\input amstex
\magnification=\magstep1
\documentstyle{amsppt}
\vsize 45pc
\NoBlackBoxes
\TagsAsMath
\define\m{\frak m}

\define\p{\frak p}

\define\OO{\Cal O}
\define\PP{\Bbb P}
\define\N{\Bbb N}
\define\Q{\Bbb Q}
\define\R{\Bbb R}
\define\Z{\Bbb Z}
\define\Spec{\operatorname{Spec }}
\define\Proj{\operatorname{Proj }}
\define\Hom{\operatorname{Hom}}
\define\Div{\operatorname{div}}
\define\der{\partial}
\define\rd#1{{\lfloor #1 \rfloor}}
\define\rup#1{{\lceil #1 \rceil}}
\define\enddemoo{$\square$ \enddemo}
\define\lra{\longrightarrow}
%\rightheadtext{F-regular and F-pure Rings}
\NoRunningHeads

%\document
\topmatter
\title F-regular and F-pure rings \\ vs. \\ log terminal and \\ 
log canonical singularities
\endtitle
\author Nobuo Hara and Kei-ichi Watanabe % $^{*}$
\endauthor
 %; Secondary ... \endsubjclass
\thanks Both authors are partially supported by Grant-in-Aid for Scientific 
Research, Japan. 
\endthanks
\address
Department of Mathematics, College of Humanities and Sciences, 
Nihon University, Setagaya-ku, Tokyo 156--0045, Japan \\
e-mail: watanabe\@math.chs.nihon-u.ac.jp %\
Department of Mathematical Sciences, Waseda University, 
3--4--1, Okubo, Shinjuku, Tokyo 169--8555, Japan
\endaddress
\endtopmatter

\document

\head Introduction \endhead

The notions of F-regular and F-pure rings are defined by Hochster and 
Huneke \cite{HH1} and Hochster and Roberts \cite{HR}, respectively, by 
using the Frobenius map in characteristic $p > 0$. These notions have some 
similarity to the notion of rational singularities defined for singularities 
of characteristic zero. But if we look more closely, we find that F-regularity 
is strictly stronger than rational singularity and that there is no implication
between F-purity and rational singularity. (There are several variants of the 
concept of ``F-regular" rings which are expected and in some cases proved to 
be equivalent to each other. In this paper, we always concern ``strongly 
F-regular" rings as in Definition 1.1 (2).) 

In dimension two, F-regular rings and F-pure rings are well investigated and 
we find very strong connection between F-regular rings and \lq\lq quotient 
singularities" \cite{W2}, \cite{Ha2}. To generalize this result to higher 
dimension, we notice that there is a notion called log terminal singularities 
for rings of characteristic zero, which is equivalent to quotient singlarities 
in dimension two. Also, the notion of log canonical singularities is defined 
similarly as log terminal singularities. Looking more closely, we find that 
F-regular (resp.\ F-pure) rings and log terminal (resp.\ log-canonical) 
singularities have very similar properties. 

The notions of F-regular and F-pure rings are also defined for rings of 
``characteristic zero."  Namely, we say that a ring essentially of finite 
type over a field of characteristic zero is of F-regular (resp.\ F-pure) 
type if its reduction to characteristic $p$ is F-regular (resp.\ F-pure) 
for infinitely many $p$. 

In this paper, we will show that a $\Q$-Gorenstein ring of F-regular (resp.\ 
F-pure) type has log terminal (resp.\ log canonical) singularities. Actually, 
we have the following result in characteristic $p > 0$, as a special case of 
our main theorem (Theorem 3.3). 
 
\proclaim{Theorem}
Let $(A,\m)$ be a $\Q$-Gorenstein normal local ring of characteristic $p > 0$ 
and $f \colon X \to Y = \Spec A$ be a proper birational morphism with $X$ 
normal. Let $E= \bigcup_{i=1}^sE_i$ be the exceptional divisor of $f$ and let 
$$K_X = f^*K_Y + \sum _{i=1}^s a_iE_i$$ 
as in $(3.1.1)$. 
If $A$ is F-pure $($resp.\ strongly F-regular$)$, then $a_i \ge -1$ 
$($resp.\ $a_i > -1)$ for every $i$. 
\endproclaim

If $A$ is a normal local ring essentially of finite type over a field of 
characteristic zero and if $f \colon X \to Y = \Spec A$ is a resolution 
of singularities of $A$, then we can apply the above theorem to ``reduction 
modulo $p$" of this map. Also, this result is applicable to strongly F-regular 
(resp.\ normal F-pure) rings of fixed characteristic $p$. It is crucially used 
in \cite{Ha2} to classify two-dimensional F-regular and normal F-pure rings 
in every characteristic $p$. 

The technique we use in our proof is that of Frobenius splitting used by 
\cite{MR}. If a ring has some splitting of Frobenius, we can lift it to an 
open set of its resolution. Then the splitting affects ``discrepancy" of 
the exceptional divisors. 

The converse of the ``F-regular" part of the above theorem has been 
established as well \cite{Ha3}, and we have a Frobenius characterization 
of log terminal singularities. Namely, a ring in characteristic zero has 
log terminal singularities if and only if it is of F-regular type and 
$\Q$-Gorenstein. See also \cite{MS2}, \cite{S1,3}. 

The other new ingredient of this paper is an attempt to generalize the notions 
of F-regular and F-pure rings to those for ``pairs." By a {\it pair}, we mean 
a pair $(A,\Delta)$ of a normal ring $A$ and a $\Q$-divisor $\Delta$ on 
$\Spec A$. The notions of log terminal and log canonical singularities are 
defined not only for normal rings (in characteristic zero) but also for pairs 
(cf.\ \cite{KMM}), and these ``singularities of pairs" play a very important 
role in birational algebraic geometry. (To wit, see Koll\'ar's lecture note 
\cite{Ko}.) Therefore we are tempted to define corresponding ``F-singularities 
of pairs." 

In this paper, we define a few variants of ``F-singularities of pairs," 
namely, strong F-regularity, divisorial F-regularity and F-purity of pairs, 
which are expected to correspond singularities of pairs called Kawamata log 
terminal (KLT), purely log terminal (PLT) and log canonical (LC). 

The significance of our main theorem (Theorem 3.3) is enhanced by considering 
F-singularities of pairs. We also prove some results on F-singularities of 
pairs, which are analoguous to the results proved for singularities of pairs 
in characteristic zero (cf.~ \cite{Ko}). If we keep in mind the (expected) 
correspondence of F-singularities and zero characteristic singularities of 
pairs, we find that Theorem 4.8 is parallel to the behavior of singularities 
of pairs under finite covering. Also, Theorem 4.9 is an analog of the 
so-called ``inversion of adjunction," and Proposition 2.10 corresponds 
to some fundamental properties of ``log canonical thresholds." 

Our main theorem for the no boundary case ($\Delta = 0$) was first proved 
by the second-named author several years ago \cite{W3}, and this result 
remained unpublished. Then the first-named author proposed a generalization 
to F-singularities of pairs, and the present paper was worked out. 

In this paper, all rings will be commutative and Noetherian. Also, since 
we are  interested only in normal local rings, {\it we assume all rings are 
integral domains}. 

\head 1. Preliminaries \endhead

Let $A$ be an integral domain of characteristic $p > 0$ and let $F \colon A 
\to A$ be the Frobenius endomorphism given by $F(a) = a^p$. We always use the 
letter $q$ for a power $q = p^e$ of $p$. 
Since $A$ is assumed to be reduced, we can identify the following three maps: 
$$
F^e \colon A  \to A, \quad 
A^q \hookrightarrow A, \quad 
A \hookrightarrow A^{1/q}. 
$$
The ring $A$ is called {\it F-finite}, if $F\colon A \to A$ (or $A 
\hookrightarrow A^{1/p}$) is a finite map. 
For example, $A$ is F-finite if it is essentially of finite type over a 
perfect field or it is a complete local ring with perfect residue field. 
{\it We always assume $A$ is F-finite} throughout this paper. 

\demo{Definition 1.1}
(1) (Hochster and Roberts \cite{HR}) $A$ is {\it F-pure} if the map 
$A \hookrightarrow A^{1/p}$ (hence $A \hookrightarrow A^{1/q}$ for every 
$q = p^e$, $e \ge 1$) splits as an $A$-module homomorphism. 

(2)\; (Hochster and Huneke \cite{HH2}) $A$ is {\it strongly F-regular} if 
for every $c \in A$ which is not in any minimal prime of $A$ (or equivalently, 
$c \ne 0$ in our case), there exists $q = p^e$ such that the $A$-module 
homomorphism $A \to A^{1/q}$ sending 1 to $c^{1/q}$ splits. 
\enddemo

F-pure and (strongly) F-regular rings have many nice properties. But here, 
we will only list the following. Interested readers could refer \cite{HH1}, 
\cite{HH2}, \cite{HR}, \cite{FW} or \cite{W2} (in the last reference, \lq\lq
F-regular" should read ``strongly F-regular").  

\remark{Remark 1.2}
(1) Let $(A,\m)$ be a local ring and $E = E_A(A/\m)$ be the injective envelope 
of the $A$-module $A/\m$. Then the splitting of the $A$-homomorphism $A \to 
A^{1/q}$ sending 1 to $c^{1/q}$ is equivalent to the injectivity of the map 
$E \to E \otimes_A A^{1/q}$ sending $\xi \in E$ to $\xi \otimes c^{1/q} \in 
E \otimes_A A^{1/q}$ (cf. Proposition 2.4). 

(2) There are notions of F-regular and weakly F-regular rings also defined 
by Hochster and Huneke \cite{HH1}, which we do not define here. On the other 
hand, a local ring $(A,\m)$ of dimension $d$ is called F-rational \cite{FW}, 
if $A$ is Cohen--Macaulay and if for every $c \ne 0 \in A$, there exists $q 
= p^e$ such that the map $cF^e \colon H_{\m}^d(A) \to H_{\m}^d(A)$, or 
equivalently, the mapping 
$$H^d_{\m}(A) = H^d_{\m}(A) \otimes_A A \to 
                   H^d_{\m}(A)\otimes_A A^{1/q} \cong H^d_{\frak m}(A^{1/q}),$$
sending $\xi \in H^d_{\m}(A)$ to $\xi \otimes c^{1/q} = (cF^e(\xi))^{1/q}$, 
is injective. 

Then the following implications hold: 
$$\text{strongly  F-regular } \Rightarrow \; \text{F-regular }
\Rightarrow \; \text{weakly  F-regular } \Rightarrow \; \text{F-rational }
\Rightarrow \; \text{normal.}
$$ 
Also, for a Gorenstein local ring, strongly F-regular is equivalent to 
F-rational. 

(3) (\cite{Mc}, \cite{Wi}) For $\Q$-Gorenstein local rings, strongly F-regular 
and weakly F-regular are equivalent. (A normal local ring is $\Q$-Gorenstein 
if its canonical module has a finite order in the divisor class group.) In 
fact, this remains true if the ring has an isolated non-$\Q$-Gorenstein point. 

(4) (\cite{HR}, \cite{HH1}, \cite{HH2}) If a ring is (strongly or weakly) 
F-regular (resp.\ F-pure), so are its pure subrings. (A ring homomorphism 
$A \to B$ is pure if for every $A$-module $M$, the map  $M \to M \otimes_A B$ 
is injective). 

(5) (\cite{MS1}, \cite{W2}) Conversely, if $(A,\m) \to (B,\frak{n})$ is 
a finite local homomorphism of normal local rings which is \'etale in 
codimension 1 and if $A$ is strongly F-regular (resp.\ F-pure), so is $B$. 
\endremark

\demo{Notation 1.3}
Let $A$ be a normal domain with quotient field $L$. A $\Q$-Weil divisor on 
$Y = \Spec A$ is a linear combination $D = \sum_{i=1}^r \alpha_iD_i$ of 
irreducible reduced subschemes $D_i \subset Y$ of codimension 1 with 
coefficients $\alpha_i \in \Q$. The round-down and round-up of $D$ 
is defined by $\rd{D} = \sum_{i=1}^r\rd{\alpha_i}D_i$ and $\rup{D} = 
\sum_{i=1}^r \rup{\alpha_i}D_i$, respectively. We also denote 
$$A(D) = \{f \in L \; \vert \; \text{div}_Y(f) + D \ge 0\}.$$ 
Clearly $A(D) = A(\rd{D})$, and this is a divisorial (i.e., finitely generated 
reflexive) submodule of $L$. Conversely, any divisorial submodule $I$ of $L$ 
is written as $I = A(D)$ for some unique {\it integer coefficient} Weil 
divisor $D$. We denote by $I^{(m)}$ the reflexive hull of $I^m$. If $I = 
A(D) = A(\rd{D})$, then $I^{(m)} = A(m\rd{D})$. 

For any ideal $I \subset A$ and $q = p^e$, we denote by $I^{[q]}$ the ideal 
of $A$ generated by the $q$th powers of elements of $I$. Also, the notation 
$(\; )^{1/q}$ will show that the module under consideration is an $A^{1/q}
$-submodule of $L^{1/q}$. For example, $(I^{[q]})^{1/q} = I \cdot A^{1/q}$. 
\enddemo

\demo{1.4. Frobenius map and its splitting}
Let $A$ be a normal domain of characteristic $p > 0$ and let $D$ be an 
effective Weil divisor on $Y = \Spec A$. Then for any $q = p^e$, we have 
a natural inclusion map $\imath \colon A \hookrightarrow A(D)^{1/q}$, which 
is identified with the $e$-times Frobenius $F^e\colon A \to A$ followed by 
the inclusion map $A \hookrightarrow A(D)$. If $D \ge qD_0$ for some 
effective divisor $D_0$, then $\imath$ is factorized as 
$$A \hookrightarrow A(D_0) = A(D_0) \otimes_A A \to 
    A(D_0) \otimes_A A^{1/q} \to A(qD_0)^{1/q} \hookrightarrow A(D)^{1/q}.$$ 
Hence, if $\imath$ splits as an $A$-module homomorphism, then $A\hookrightarrow
A(D_0)$ also splits, and this implies $D_0 = 0$. This kind of argument is used 
frequently in the sequel. Note also that a splitting is preserved under 
localization, in particular, a localization of a strongly F-regular ring 
is again strongly F-regular. 
\enddemo

\head 2. F-regularity and F-purity of Pairs \endhead

\demo{Definition 2.1}
Let $A$ be an F-finite normal domain of characteristic $p > 0$ and 
$\Delta$ be an effective $\Q$-Weil divisor on $Y = \Spec A$. 

(1) We say that the pair $(A,\Delta)$ is {\it F-pure} if the inclusion map 
$A \hookrightarrow A((q-1)\Delta)^{1/q}$ splits as an $A$-module homomorphism 
for every $q = p^e$. 

(2) $(A,\Delta)$ is {\it strongly F-regular} if for every nonzero element 
$c \in A$, there exists $q = p^e$ such that $c^{1/q}A \hookrightarrow A((q-1)
\Delta)^{1/q}$ splits as an $A$-module homomorphism. 

(3) $(A,\Delta)$ is {\it divisorially F-regular} if for every nonzero element 
$c\in A$ which is not in any minimal prime ideal of $A(-\rd{\Delta}) \subseteq 
A$, there exists $q = p^e$ such that $c^{1/q}A \hookrightarrow A((q-1)\Delta
)^{1/q}$ splits as an $A$-module homomorphism. 
\enddemo

Definition 2.1 includes Definition 1.1 as the special case $\Delta = 0$. 
Namely, $A$ is F-pure if and only if $(A,0)$ is F-pure, and $A$ is strongly 
F-regular if and only if $(A,0)$ is strongly F-regular, or equivalently, 
$(A,0)$ is divisorially F-regular. 

We collect some basic properties in the following. 

\proclaim{Proposition 2.2}
Let $A$ be an F-finite normal domain of characteristic $p > 0$ and $\Delta$ 
be an effective $\Q$-Weil divisor on $Y = \Spec A$. 
\roster
\item For a pair $(A,\Delta)$ as above, the following implications hold. 
$$\text{
strongly F-regular $\Rightarrow$ divisorially F-regular $\Rightarrow$ F-pure. 
}$$
\item $(A,\Delta)$ is strongly F-regular if and only if for every nonzero 
element $c\in A$, there exists $q'$ such that the map $c^{1/q}A \hookrightarrow
A(q\Delta)^{1/q}$ splits as an $A$-module homomorphism for all $q= p^e\ge q'$. 
\item $(A,\Delta)$ is strongly F-regular if and only if $(A,\Delta)$ is 
divisorially F-regular and $\rd{\Delta} = 0$. 
\item If $(A,\Delta)$ is  F-pure, then $\rup{\Delta}$ is reduced, i.e., the 
coefficient of $\rup{\Delta}$ in every irreducible component is equal to $1$. 
\item If $(A,\Delta)$ is F-pure $($resp.\ strongly or divisorially 
F-regular$)$, so is $(A,\Delta')$ for every effective $\Q$-Weil divisor 
$\Delta' \le \Delta$.
\endroster
\endproclaim

\demo{Proof}
(1) It is obvious that strongly F-regular implies divisorially F-regular. 
Also, divisorially F-regular is equivalent to F-pure if $\dim A = 1$. 
If $\dim A \ge 2$, we can choose a non-unit $c \ne 0$ of $A$ which is 
not in any minimal prime of $A(-\rd{\Delta})$. Then the map $(c^q)^{1/q'}A 
\hookrightarrow A((q'-1)\Delta)^{1/q'}$ does not split for every $q' \le q$, 
because it factors through $cA \hookrightarrow A$ (cf.~1.4). Hence, if $(A,
\Delta)$ is divisorially F-regular and if $q$ is any power of $p$, then 
there is a power $q' > q$ such that $(c^q)^{1/q'}A \hookrightarrow 
A((q'-1)\Delta)^{1/q'}$ splits. This implies that $A \hookrightarrow 
A((q'-1)\Delta)^{1/q'}$ splits, and so does $A \hookrightarrow A((q-1)
\Delta)^{1/q}$. Consequently, $(A,\Delta)$ is F-pure. 

(2) The sufficiency is clear. To show the necessity, let $c \in A$ be any 
nonzero element, and choose $d \ne 0 \in A(-\Delta)$ so that $dA(\rup{q\Delta})
\subset A((q-1)\Delta)$ for every $q = p^e$. If $(A,\Delta)$ is strongly 
F-regular, then there exists a power $q'$ of $p$ such that $A @>(cd)^{1/q'}>> 
A((q'-1)\Delta)^{1/q'}$ splits. Since this map is factorized into $A \overset
{c^{1/q'}}\to\lra A(\rup{q'\Delta})^{1/q'} \overset{d^{1/q'}}\to\lra A((q'-1)
\Delta)^{1/q'}$, the map $A \overset{c^{1/q'}}\to\lra A(\rup{q'\Delta})^{1/q'}$
also splits. On the other hand, since $A$ is F-pure, the map $A^{1/q'} 
\hookrightarrow A^{1/qq'}$ splits for all $q = p^e$, and so does $A(\rup{q'
\Delta})^{1/q'} \hookrightarrow A(q\rup{q'\Delta})^{1/qq'}$, too. 
Hence the map $A \overset{c^{1/q'}}\to\lra A(q\rup{q'\Delta})^{1/qq'}$ 
splits for all $q = p^e$. Since this map is factorized into $A @>c^{1/qq'}>> 
A(qq'\Delta)^{1/qq'} @>c^{(q-1)/qq'}>> A(q\rup{q'\Delta})^{1/qq'}$, the map 
$A @>c^{1/q}>> A(q\Delta)^{1/q}$ splits for all $q = p^e \ge q'$. 

To prove (3) (resp.\ (4)), assume to the contrary that $\Delta$ has a component
$\Delta_0$ with coefficient $\ge 1$ (resp.\ $> 1$). Then there is a $q = p^e$ 
such that the coefficient of $q\Delta$ (resp.\ $(q-1)\Delta$) in $\Delta_0$ 
is at least $q$. Then the map $A \hookrightarrow A(q\Delta)^{1/q}$ (resp.\ $A 
\hookrightarrow A((q-1)\Delta)^{1/q}$) factors through $A \hookrightarrow A(
\Delta_0) \hookrightarrow A(q\Delta_0)^{1/q}$, which does not split. This 
implies that $(A,\Delta)$ cannot be strongly F-regular (resp.\ F-pure). 

(5) The only being non-trivial is the assertion for divisorial F-regularity 
in the case $\text{Supp}(\rd{\Delta}) \setminus \text{Supp}(\rd{\Delta'}) 
\ne \emptyset$. To prove this, we may assume without loss of generality that 
there is a unique irreducible component $\Delta_0$ of $\rd{\Delta}$ such that 
$\Delta_0 \not\subseteq \text{Supp}(\rd{\Delta'})$. Let $c \ne 0 \in A$ be 
any element which is in $\p = A(-\Delta_0)$ but not in any minimal prime of 
$A(-\rd{\Delta'})$, and let $\nu = v_{\p}(c)$, the value of $c$ at $\p$. Then 
by prime avoidance, we can choose $d\ne 0\in cA(\nu\Delta_0)$ which is not in 
any minimal prime of $A(-\rd{\Delta})$. This implies that $dA((q-1)\Delta') 
\subseteq cA((q-1)\Delta)$ for $q \gg 0$. Now, if $(A,\Delta)$ is divisorially 
F-regular, then there exists $q = p^e$ such that $A \overset{d^{1/q}}\to\lra 
A((q-1)\Delta)^{1/q}$ splits, and we can show as in (1) that this is true for 
infinitely many $q$. Since this map factors into $A \overset{c^{1/q}}\to\lra 
A((q-1)\Delta')^{1/q} \overset{(d/c)^{1/q}}\to\lra A((q-1)\Delta)^{1/q}$ if 
$q \gg 0$, the map $A \overset{c^{1/q}}\to\lra A((q-1)\Delta')^{1/q}$ splits 
for some $q$. 
\enddemoo

\remark{Remark 2.3}
In the definition of ``F-pure," it seems apparently natural to refer the map 
$A \hookrightarrow A(q\Delta)^{1/q}$ instead of $A \hookrightarrow A((q-1)
\Delta)^{1/q}$. (We have seen in Proposition 2.2 (2) that this makes 
no difference for strong F-regularity.) But it does make crucial 
difference for F-purity. Let us say $(A,\Delta)$ is {\it strongly 
F-pure} if $A \hookrightarrow A(q\Delta)^{1/q}$ splits for every $q = p^e$. 
Then the proof of Proposition 2.2 (3) shows that if $(A,\Delta)$ is strongly 
F-pure, then $\rd{\Delta} = 0$. But this is too much stronger than what we 
want for ``F-purity." Note also that we have the implication ``strongly 
F-regular $\Rightarrow$ strongly F-pure $\Rightarrow$ F-pure," and that 
there is no implication between strong F-purity and divisorial F-regularity. 
\endremark

\proclaim{Proposition 2.4}
Let $(A,\m)$ be a $d$-dimensional F-finite normal local ring of characteristic 
$p > 0$ and let $\Delta$ be an effective $\Q$-Weil divisor on $Y = \Spec A$. 
Then: 
\roster
\item $(A,\Delta)$ is F-pure if and only for every $q = p^e$, the induced 
$e$-times Frobenius map $F^e \colon H_{\m}^d(K_A) \to H_{\m}^d(A(qK_A+(q-1)
\Delta))$ is injective. 
\item $(A,\Delta)$ is strongly F-regular if and only if for every $c \ne 0 
\in A$, there exists $q = p^e$ such that $cF^e \colon H_{\m}^d(K_A) \to 
H_{\m}^d(A(qK_A+(q-1)\Delta))$ is injective. 
\item $(A,\Delta)$ is divisorially F-regular if and only if for every $c \ne 
0 \in A$ which is not in any minimal prime of $A(-\rd{\Delta})$, there exists 
$q=p^e$ such that $cF^e\colon H_{\m}^d(K_A) \to H_{\m}^d(A(qK_A+(q-1)\Delta))$ 
is injective. 
\endroster
\endproclaim

\demo{Proof}
The map $c^{1/q}A \hookrightarrow A((q-1)\Delta)^{1/q}$ splits as an 
$A$-module homomorphism if and only if the map 
$$\Hom_A(A((q-1)\Delta)^{1/q},A) \overset{c^{1/q}}\to\lra 
                        \Hom_A(A((q-1)\Delta)^{1/q},A) \to \Hom_A(A,A) = A$$ 
is surjective. By the local duality, the Matlis dual of $\Hom_A(A((q-1)\Delta
)^{1/q},A) \cong \Hom_A(A(qK_A+(q-1)\Delta)^{1/q},K_A)$ is $H_{\m}^d(A(qK_A+
(q-1)\Delta)^{1/q})$, so that the surjectivity of the above map is equivalent 
to the injectivity of 
$$H_{\m}^d(K_A) \to H_{\m}^d(A(qK_A+(q-1)\Delta)^{1/q}) 
             \overset{c^{1/q}}\to\lra H_{\m}^d(A(qK_A+(q-1)\Delta)^{1/q}).$$ 
Since this map is identified with $cF^e \colon H_{\m}^d(K_A) \to H_{\m}^d(A
(qK_A+(q-1)\Delta))$, we obtain the assertions (1), (2) and (3) at once. 
\enddemoo

\proclaim{Corollary 2.5}
If $A$ is a normal toric ring over a perfect field $k$ of characteristic 
$p > 0$ and if $\Delta$ is a reduced toric divisor, then the pair $(A,\Delta)$ 
is F-pure. 
\endproclaim

\demo{Proof}
Let $M = \Z^d$, $N = \Hom_{\Z}(M,\Z)$ and denote the duality pairing of 
$M_{\R} = M \otimes_{\Z} \R$ with $N_{\R} = N \otimes_{\Z} \R$ by $\langle \;
,\; \rangle \colon M_{\R} \times N_{\R} \to \R$. Let $A$ be the toric ring 
defined by a rational polyhedral cone $\sigma \subset N_{\R}$ and let $D_1,
\dots ,D_s$ be the toric divisors of $\Spec A$ corresponding to the primitive 
generators $n_1, \dots , n_s \in N$ of $\sigma$, respectively. To show the 
F-purity of the pair $(A,\Delta)$, we may assume that $\Delta = \sum_{i=1}^s 
D_i$, by Proposition 2.2 (5). Then one has $K_A = -\Delta$, so that $A(qK_A + 
(q-1)\Delta) = K_A$. Therefore 
$$H_{\m}^d(K_A) = H_{\m}^d(A(qK_A + (q-1)\Delta)) 
                = \bigoplus_{m \in -\check{\sigma} \cap M} k \cdot x^m,$$ 
where $-\check{\sigma} = \{m \in M_{\R} \; \vert \; \langle m,n_i \rangle \le 
0 \text{ for } i = 1,\dots ,s\}$, and the $e$-times Frobenius map 
$F^e \colon H_{\m}^d(K_A) \to H_{\m}^d(A(qK_A+(q-1)\Delta))$ is given by 
$F^e(x^m) = x^{qm}$ $(m \in -\check{\sigma} \cap M)$. Since this map is 
injective, $(A,\Delta)$ is F-pure by Proposition 2.4 (1). 
\enddemoo

\proclaim{Proposition 2.6} {\rm (Fedder-type criteria \cite{F})}
Let $(R,\m)$ be an F-finite regular local ring of characteristic $p > 0$ and 
$I \subset R$ be an ideal such that $A= R/I$ is normal. Let $g \in R \setminus 
I$ be an element whose image $\bar{g} \in A$ defines a reduced divisor $\Div_Y
(\bar{g})$ on $Y = \Spec A$, and consider a $\Q$-divisor $\Delta = t \cdot 
\Div_Y(\bar{g})$ for nonnegative $\alpha\in \Q$. We put $r_e = \rd{t(q-1)}$ 
for each $q = p^e$. Then: 
\roster
\item $(A,\Delta)$ is F-pure if and only if $g^{r_e}(I^{[q]}:I) \not\subseteq 
\m^{[q]}$ for all $q = p^e$. 
\item $(A,\Delta)$ is strongly F-regular if and only if for every $c \in R 
\setminus I$, there exists $q=p^e$ such that $cg^{r_e}(I^{[q]}:I)\not\subseteq 
\m^{[q]}$. 
\item $(A,\Div_Y(\bar{g}))$ is divisorially F-regular if and only if for every 
$c \in R \setminus I$ which is not in any minimal prime of $I + gR$, there 
exists $q = p^e$ such that $cg^{q-1}(I^{[q]}:I) \not\subseteq \m^{[q]}$. 
\endroster
\endproclaim

\demo{Proof}
The proof is essentially the same as those in \cite{F}, \cite{Gl}. First, 
we may assume without loss of generality that $(R,\m)$ is a {\it complete} 
regular local ring. We will describe the map $cF^e \colon H_{\m}^d(K_A) \to 
H_{\m}^d(A(qK_A+(q-1)\Delta))$ in Proposition 2.4. Since $A((q-1)\Delta) = 
\bar{g}^{-r_e}A$ and since $H_{\m}^d(K_A^{(q)}) \cong E_A(A/\m A) \otimes_A 
A^{1/q}$ by \cite{W2, Theorem 2.5}, This map is viewed as 
$$c\bar{g}^{r_e}F^e \colon E_A(A/\m A) \to E_A(A/\m A) \otimes_A A^{1/q}.
\tag{2.6.1}$$ 

Let $E_R = E_R(R/\m)$, $E_A = E_A(A/\m A)$ and let $n = \dim R$. Then $E_R 
\cong H_{\m}^n(R)$ since $R$ is Gorenstein, and we can identify $E_R$ with 
$E_R \otimes_R R^{1/q} \cong H_{\m}^n(R^{1/q})$ via the identification of 
$R$ with $R^{1/q}$. Also, $E_A \cong \Hom_R(R/I,E_R) \cong (0:I)_{E_R} 
\subset E_R$, and $E_A$ and $E_R$ have the 1-dimensional socle in common. 
Let $z$ be a generator of the socle of $E_R$ and let $z'$ be the corresponding 
socle generator in $E_A$. 

Now, $R \hookrightarrow R^{1/q}$ is flat since $R$ is regular, so that via 
the identification $R^{1/q} \cong R$, we have 
$$\aligned
E_A \otimes_R R^{1/q} & 
\cong \Hom_{R^{1/q}}(R/I \otimes_R R^{1/q},E_R \otimes_R R^{1/q}) \\ & 
\cong \Hom_R(R/I^{[q]},E_R) \cong (0:I^{[q]})_{E_R} 
\endaligned \tag{2.6.2}$$ 
in $E_R \otimes_R R^{1/q} \cong E_R$. Accordingly we have 
$$E_A \otimes_A A^{1/q} \cong E_A \otimes_R R^{1/q} \otimes_{R^{1/q}} A^{1/q}
                        \cong (0:I^{[q]})_{E_R} \otimes_R A 
                        \cong \frac{(0:I^{[q]})_{E_R}}{I(0:I^{[q]})_{E_R}},$$ 
and the image of $z' \in E_A$ by the map (2.6.1) is $cg^{r_e}z^q$ mod 
$I(0:I^{[q]})_{E_R}$, where $z^q$ denotes the image of $z$ by the $e$-times 
Frobenius on $E_R \cong H_{\m}^n(R)$. We also know from (2.6.2) (by putting 
$I = \m$) that $z^q \in E_R$ generates $(0:\m^{[q]})_{E_R}$. 

Now the map (2.6.1) is injective if and only if $cg^{r_e}z^q \notin I(0:I^{[q]}
)_{E_R}$ if and only if $cg^{r_e}(0:\m^{[q]})_{E_R} \not\subseteq I(0:I^{[q]}
)_{E_R}$. Since $I(0:I^{[q]})_{E_R} = (0:(I^{[q]}:I)_R)_{E_R}$ by the Matlis 
duality, this is equivalent to saying that $cg^{r_e}(I^{[q]}:I) \not\subseteq 
\m^{[q]}$. 
\enddemoo

\proclaim{Corollary 2.7} {\rm (cf.\ \cite{Ha1})}
Let $(R,\m)$ be an F-finite regular local ring of characteristic $p > 0$, and 
let $f_1, \dots ,f_s, g \in \m$ be an $R$-regular sequence. Let $A = R/(f_1,
\dots ,f_s)$ and consider an effective $\Q$-divisor $\Delta = t\cdot\Div_Y
(\bar{g})$ on $Y = \Spec A$. Put $r_e = \rd{t(p^e-1)}$. Then: 
\roster
\item $(A,\Delta)$ is F-pure if and only if $(f_1 \cdots f_s)^{q-1}g^{r_e} 
\notin \m^{[q]}$ for all $q = p^e$. 
\item $(A,\Delta)$ is strongly F-regular if and only if 
$$\bigcap_{q=p^e, e\in \N} \m^{[q]}:(f_1 \cdots f_s)^{q-1}g^{r_e} 
                                   = (f_1,\dots ,f_s) \quad \text{in }R.$$ 
\item $(A,\Div_Y(\bar{g}))$ is F-pure $($resp.\ divisorially F-regular$)$ if 
and only if $A/\bar{g}A = R/(f_1,\dots ,f_s,g)$ is F-pure $($resp.\ strongly 
F-regular$)$. 
\endroster
\endproclaim

\demo{Proof}
The proof is easy, but we remark one point which might be overlooked. 
Put $D = \Div_Y(\bar{g})$. Then the condition $(f_1 \cdots f_s)^{q-1}g^{r_e} 
\notin \m^{[q]}$ holds if and only if $A\hookrightarrow A(\rd{t(q-1)}D)^{1/q}$ 
splits as an $A$-module homomorphism by the preceding argument. So, to prove 
(1), we have to show the equivalence of the following two conditions. 
\roster
\item"{(i)}" $A \hookrightarrow A((q-1)\Delta)^{1/q} = A(\rd{(q-1)tD})^{1/q}$ 
splits for all $q=p^e$. 
\item"{(ii)}" $A \hookrightarrow A(\rd{t(q-1)}D)^{1/q}$ splits for all $q=p^e$.
\endroster
The implication (i) $\Rightarrow$ (ii) is clear since $\rd{(q-1)tD} \ge 
\rd{(q-1)t}D$. But (ii) $\Rightarrow$ (i) is not apparently clear if $D$ is 
non-reduced. To show this implication, let $q$ be any power of $p$, and choose 
$q'$ such that $q'\rd{(q-1)tD} \le \rd{t(qq'-1)}D$. Then by (ii), the map $A 
\hookrightarrow A(\rd{t(qq'-1)}D)^{1/qq'}$ splits, and this map is factorized 
into $A \hookrightarrow A(\rd{(q-1)tD})^{1/q} \hookrightarrow A(q'\rd{(q-1)tD}
)^{1/qq'} \hookrightarrow A(\rd{t(qq'-1)}D)^{1/qq'}$. Hence $A \hookrightarrow 
A(\rd{(q-1)tD})^{1/q}$ splits. 
\enddemoo

\remark{Remark 2.8}
We have Fedder-type criteria also for strong F-purity (see Remark 2.3 for 
a definition), by setting $r_e = \rd{tp^e}$ in (1) of Proposition 2.6 and 
Corollary 2.7. The criteria for strong F-regularity also works if we put 
$r_e = \rd{tp^e}$. 
\endremark

\demo{Example 2.9}
(1) Let $A$ be a regular local ring and let $\Delta = t \cdot \Div(x_1 \cdots 
x_i)$ on $\Spec A$, where $x_1, \dots ,x_i$ are part of regular parameters of 
$A$. If $t \le 1$ (resp.\ $t < 1$), then $(A,\Delta)$ is F-pure (resp.\ 
strongly F-regular). 

(2) Let $A = k[[X,Y,Z]]/(XY-Z^2)$ and denote the images of $Z$ in $A$ by $z$. 
Let $\Delta = t \cdot \Div(z)$. Corollary 2.7 tells us that if $t \le 1$ 
(resp.\ $t < 1$), then $(A,\Delta)$ is F-pure (resp.\ strongly F-regular). 

(3) Let $A = k[[x,y]]$ and let $\Delta = \frac56 \Div(x^2-y^3)$. Then in any 
characteristic $p > 0$, $(A,\Delta)$ is F-pure but not strongly F-regular. 
But $(A,\Delta)$ is strongly F-pure if and only if $p \equiv 1$ mod $3$. 
\enddemo

The following is a variant of Fedder's result \cite{F}, see also \cite{Ko, 
Lemma 8.10}. 

\proclaim{Proposition 2.10}
Let $A = k[[x_1, \dots ,x_d]]$ be a $d$-dimensional complete regular local 
ring over a perfect field $k$ of characteristic $p > 0$ and let $f \in A$ be a 
nonzero element of multiplicity $n$, i.e., $f \in \m^{n} \setminus \m^{n+1}$, 
where $\m$ is the maximal ideal of $A$. For nonnegative $t \in \Q$ we have: 
\roster
\item If $t \le 1/n$, then $(A, t\cdot\Div (f))$ is F-pure. 
\item If $(A,t\cdot\Div (f))$ is F-pure $($resp.\ strongly F-regular$)$, 
then $t \le d/n)$ $($resp.\ $t < d/n)$. 
\item Assume that $t < \operatorname{min} \{1,d/n\}$ and that the initial 
term $f_n$ of $f$ defines a smooth subvariety of $\PP^{d-1}$. 
Then in characteristic $p \gg 0$, $(A,t\cdot\Div(f))$ is strongly F-pure. 
More precisely, if the ideal $J= (\der f_n/\der x_1,\dots ,\der f_n/\der x_d)$ 
contains $(x_1^{m_1}, \dots ,x_d^{m_d})$ in $A$, then $(A,t\cdot\Div(f))$ 
is strongly F-pure, or otherwise $p < (m_1+ \cdots +m_d)/(d-nt)$. 
\endroster
\endproclaim

\demo{Remark} Note that if $n = \deg f_n$ is not divisible by $p$, then 
$f_n \in J$, so that $J$ is equal to the Jacobian ideal $(f_n,J)$ of $V = 
(f_n = 0) \subset \PP^{d-1}$. So the smoothness of $V$ implies that 
$(x_1^{m_1}, \dots ,x_d^{m_d}) \subseteq J$ for some $m_1, \dots ,m_d \in \N$. 
\enddemo

\demo{Proof}
The assertions (1) and (2) immediately follows from Corollary 2.7. To prove 
(3) we use the Fedder-type critrion for strong F-purity (2.8): 
$$(A,t\cdot\Div (f)) \text{ is strongly F-pure } \Longleftrightarrow 
                  \; f^{r_e} \notin \m^{[p^e]} \text{ for all } e \in \N,$$ 
where $r_e = \rd{tp^e}$. Since $f_n^{r_e} \notin \m^{[p^e]}$ implies 
$f^{r_e} \notin \m^{[p^e]}$, we may assume that $f$ is a homogeneous polynomial
of degree $n$ in $x_1, \dots ,x_d$. 

Assume $p \ge \mu := (m_1+ \cdots +m_d)/(d-nt)$ and let $j_e$ be the integer 
such that $f^{j_e+1} \in \m^{[p^e]}$ but $f^{j_e} \notin \m^{[p^e]}$. Then it 
is sufficient to prove that $r_e \le j_e$ for every $e \ge 0$. Assume to the 
contrary that there exists an $e$ such that $r_e>j_e$, and choose the smallest 
one among all such $e$. Then $e> 0$, since $r_0= 0$ by the assumption $t < 1$. 
Also, by the minimality of $e$, we have $f^{r_{e-1}} \notin \m^{[p^{e-1}]}$, 
and this implies that $f^{pr_{e-1}}\notin \m^{[p^e]}$. Indeed, if $f^{pr_{e-1}}
\in \m^{[p^e]}$, then $f^{r_{e-1}} = (f^{pr_{e-1}})^{1/p} \in \m^{[p^{e-1}]}
A^{1/p} \cap A = \m^{[p^{e-1}]}$, because $A\to A^{1/p}$ is pure. Hence 
$pr_{e-1} \le j_e < r_e$, and $j_e+1$ is not divisible by $p$ since $r_e - 
pr_{e-1} \le p-1$. Therefore, by differentiating $f^{j_e+1} \in \m^{[p^e]}$ 
by $x_i$, we have $f^{j_e}\der f/\der x_i \in \m^{[p^e]}$ for $i = 1, \dots ,
d$, so that 
$$f^{j_e} \in \m^{[p^e]}:J \subseteq \m^{[p^e]}:(x_1^{m_1}, \dots , x_d^{m_d}) 
               \subseteq \m^{[p^e]} + (x_1^{p^e-m_1} \cdots x_d^{p^e-m_d})A.$$ 
Then $f^{j_e}$ must have a nonzero term in $(x_1^{p^e-m_1}\cdots x_d^{p^e-m_d})
A$, since $f^{j_e} \notin \m^{[p^e]}$. Hence $p^e d - \sum_{i=1}^d m_i = \deg 
f^{j_e}= nj_e< nr_e \le ntp^e$. This inequality, together with the assumption 
$t < d/n$, implies $p^e < \mu$, which contradicts to $p \ge \mu$. 
\enddemoo

\head 3. Main Theorem \endhead

First, we recall the definition of log terminal and log canonical 
singularities. 

Let $f \colon X \to Y$ be a proper birational morphism of normal varieties 
over a field $k$ and let $E = \bigcup_{i=1}^s E_i$ be the exceptional divisor 
of $f$ with irreducible components $E_1, \dots ,E_s$. For a $\Q$-Weil divisor 
$D$ on $Y$ (resp.\ on $X$), we denote by $f^{-1}_*(D)$ (resp.\ $f_*D$) the 
strict transform of $D$ in $X$ (resp.\ in $Y$).  A $\Q$-Weil divisor $D$ on 
$Y$ is said to be {\it $\Q$-Cartier} if $rD$ is a Cartier divisor for some 
integer $r>0$. Then the pull-back $f^*(rD)$ of $rD$ is also a Cartier divisor 
on $X$, and we can define the pull-back $f^*D$ of $D$ as a $\Q$-Cartier 
divisor by $f^*D = \frac1r f^*(rD)$. 

We denote the dualizing sheaves of $X$ and $Y$ by $\omega_X$ and $\omega_Y$, 
respectively, and fix $\omega_X$ and $\omega_Y$ as divisorial subsheaves of 
the rational function field $k(X) = k(Y)$ so that they coincide with each 
other outside the exceptional locus of $f$. Then the canonical divisor $K_X$ 
of $X$ (resp.\ $K_Y$ of $Y$) is also fixed by $\omega_X = \OO_X(K_X)$ (resp.\ 
$\omega_Y = \OO_Y(K_Y)$), and the strict transform of $K_X$ in $Y$ is $f_*K_X 
= K_Y$. 

Now let $\Delta$ be an effective $\Q$-Weil divisor on $Y$ such that $K_Y + 
\Delta$ is $\Q$-Cartier and denote by $\tilde \Delta$ the strict transform 
$f^{-1}_*(\Delta)$ of $\Delta$ in $X$. Let $r > 0$ be an integer such that 
$r(K_Y + \Delta)$ is Cartier. Then $r(K_X + \tilde{\Delta})$ and $r \cdot f^*
(K_Y + \Delta) = f^*(r(K_Y + \Delta))$ have integer coefficients, and coincide 
with each other outside the exceptional locus of $f$. Hence $r(K_X + \tilde{
\Delta}) = r \cdot f^*(K_Y + \Delta) + \sum_{j=1}^s b_j E_j$ for some $b_1,
\dots ,b_s \in \Z$, and we have 
$$K_X + \tilde{\Delta} = f^*(K_Y + \Delta) + \sum_{j=1}^s a_j E_j, 
     \quad \text{where} \; a_j = \frac{b_j}r \; (j=1,\ldots ,s). \tag{3.1.1}$$ 
%If we denote $E \cup \text{Supp}\tilde{\Delta}$ by $F = \bigcup_{j=1}^s F_j$, 
%we can also write $$K_X = f^*(K_Y+\Delta) + \sum_{j=1}^t a_jF_j \tag{3.1.2}$$ 
%for some $a_j \in \Q$ ($j = 1,\dots ,t$). 
We call $a_j \in \Q$ the {\it discrepancy} of $E_j$ with respect to 
$(Y,\Delta)$. 
In the following definition, we keep this notation for any desingularization 
$f \colon X \to Y$ under consideration. 

\demo{Definition 3.1}
Let $Y$ be a normal variety of characteristic zero, $\Delta$ be an effective 
$\Q$-Weil divisor, and assume that $K_Y + \Delta$ is $\Q$-Cartier. 

(1) $(Y,\Delta)$ is said to be {\it Kawamata log terminal} (or {\it KLT} for 
short), if $\rd{\Delta} = 0$ and if for every resolution of singularities 
$f \colon X \to Y$ and every $f$-exceptional divisor $E_j$, the discrepancy 
$a_j$ defined in (3.1.1) satisfies $a_j > -1$. 

(2) $(Y,\Delta)$ is {\it purely log terminal} (or {\it PLT} for short), if for 
every resolution of singularities $f \colon X \to Y$ and every $f$-exceptional 
divisor $E_j$, one has $a_j > -1$. 

(3) $(Y,\Delta)$ is {\it log canoninal\/} (or {\it LC\/} for short), if for 
every resolution of singularities $f \colon X \to Y$ and every $f$-exceptional 
divisor $E_j$, one has $a_j \ge -1$. 
\enddemo

When the pair $(Y,0)$ is KLT, or equivalently, PLT in this case (resp.\ LC), 
we say that $Y$ has log terminal (resp.\ log canonical) singularities. 

\remark{Remark 3.2}
(1) Clearly, we have the implications ``KLT $\Rightarrow$ PLT $\Rightarrow$ 
LC," and if $(Y,\Delta)$ is LC, then $\rup{\Delta}$ is reduced, i.e., every 
coefficients of $\Delta$ is less than or equal to 1. 
Also, conditions (1) and (3) of Definition 3.1 is checked by referring {\it 
some} log resolution of $(Y,\Delta)$, that is, a resolution of singularities 
$f \colon X \to Y$ such that the union of the exceptional set and $\text{Supp}
(f^{-1}\Delta)$ is a simple normal crossing divisor. Namely, $(Y,\Delta)$ is 
KLT (resp.\ LC) if and only if $\rd{\Delta} = 0$ (resp.\ $\rup{\Delta}$ is 
reduced) and there exists a log resolution $f \colon X \to Y$ such that 
$a_j > -1$ (resp.\ $a_j \ge -1$) for every $f$-exceptional divisor $E_j$. 
Also, $(Y,\Delta)$ is PLT if and only if there exists a log resolution 
$f \colon X \to Y$ such that $f^{-1}_*\rd{\Delta}$ is smooth and that 
$a_j > -1$ for every $f$-exceptional divisor $E_j$.

(2) Log terminal (resp.\ log canonical) singularities behave similarly as 
strongly F-regular (resp.\ F-pure) rings under finite covers. Namely, if 
$(A,\m) \to (B,\frak n)$ is a finite local homomorphism of normal local 
rings which is \'etale in codimension 1, then $A$ is log terminal (resp.\ 
log canonical) if and only if so is $B$. 

(3) If $A$ has log terminal singularities, then $A$ has rational singularities.
Conversely, if $A$ has Gorenstein rational singularities, then $A$ has log 
terminal singularities. 
\endremark

\proclaim{Theorem 3.3}
Let $(A,\m)$ be an F-finite normal local ring of characteristic $p > 0$ and 
$\Delta$ be an effective $\Q$-Weil divisor on $Y = \Spec A$ such that $K_Y 
+ \Delta$ is $\Q$-Cartier. Let $f \colon X \to Y = \Spec A$ be a proper 
birational morphism with $X$ normal. Let $E = \bigcup _{j=1}^s E_j$ be the 
exceptional divisor of $f$ and let 
$$K_X + \tilde{\Delta} = f^*(K_Y + \Delta) + \sum _{j=1}^s a_j E_j$$ 
as in $(3.1.1)$. 
If the pair $(A,\Delta)$ is F-pure $($resp.\ divisorially F-regular$)$, 
then $a_j \ge -1$ $($resp.\ $a_j > -1)$ for every $j = 1,\dots ,s$. 
\endproclaim

We will give two different proofs to this theorem, in which we relax 
the properness assumption of $f \colon X \to Y = \Spec A$ as follows: $f$ 
factorizes into $f \colon X \hookrightarrow \bar{X} \overset\bar{f}\to\lra Y$, 
where $\bar{f} \colon \bar{X} \to Y$ is a proper birational morphism with 
$\bar X$ normal, and $X$ is an open subset of $\bar X$ with $\text{codim}
(\bar{X} \setminus X,\bar{X}) \ge 2$. This assumption is preserved if we 
replace $X$ by the nonsingular locus $X_{\text{reg}}$ of $X$, and we do not 
lose any information about the discrepancies $a_j$ by this replacement, since 
$X_{\text{reg}}$ intersects every $E_j$. Hence we may assume without loss of 
generality that {\it $X$ is Gorenstein} in what follows. 

\demo{First proof}
We will divide the proof into six steps. First (since $\bar f$ is proper and 
birational), $f$ is an isomorphism on an open set $U$ of \; $Y$ with $\text{
codim}(Y \setminus U,Y) \ge 2$. Since $K_A^{(i)}$ is reflexive, we have 
$$H^0(X,\omega_X^{\otimes i}) \subseteq 
                                   H^0(U,\omega_U^{\otimes i}) = K_A^{(i)}$$
for every $i \in \Z$. Also, we can choose a nonzero element $b \in A$ such 
that 
$$b \cdot K_A^{(-i)} \subseteq H^0(X,\omega_X^{\otimes (-i)}) \quad 
                     \text{for each } i = 0,1,\ldots ,r-1. \tag{3.3.1}$$ 

(1) Since $X$ is Gorenstein, we have 
$$\aligned
{\Cal H}om_{\OO_X}(\OO_X^{1/q},\OO_X) & \cong 
{\Cal H}om_{\OO_X}(\OO_X^{1/q},\omega_X) \otimes_{\OO_X} \omega_X^{-1} \\
 & \cong (\omega_X)^{1/q} \otimes \omega_X^{-1} \cong \OO_X((1-q)K_X)^{1/q}. 
\endaligned \tag{3.3.2}$$
by the adjunction formula for the finite map $\OO_X \hookrightarrow 
\OO_X^{1/q}$. By a similar argument on $U$, we also have an isomorphism 
$$\alpha \colon \Hom_A(A^{1/q},A) \overset{\sim}\to\lra (K_A^{(1-q)})^{1/q},$$ 
which induces 
$\Hom_A(A((q-1)\Delta)^{1/q},A) \cong A((1-q)K_A - \rd{(q-1)\Delta})^{1/q}$ 
and is compatible with (3.3.2). 

(2) Let us fix an embedding $\OO_X(iK_X) \hookrightarrow L$, where $L = k(X)$ 
is the rational function field of $X$, and put $A(r(K_A + \Delta)) = wA$. 
Then $A((mr+i)(K_A + \Delta)) = w^mA(i(K_A + \Delta))$, and 
$\omega_X^{\otimes r} = \OO_X(rK_X) = \OO_X(r(f^*(K_Y+\Delta)+\sum a_jE_j-
\tilde{\Delta})) = f^*w \cdot \OO_X(r\sum a_jE_j-\tilde{\Delta})$. 

(3) Now, assume that $(A,\Delta)$ is F-pure and let $\phi \colon A((q-1)\Delta
)^{1/q}\to A$ be a splitting of $A\hookrightarrow A((q-1)\Delta)^{1/q}$. 
Then $\phi$ induces a splitting $\tilde{\phi} \colon L^{1/q} \to L$ of 
$L \hookrightarrow L^{1/q}$. Let $X' = X \setminus Z$, where $Z = \text{Supp}
((\tilde{\phi}(\OO_X^{1/q}) + \OO_X)/\OO_X)$. Then $X'$ is an open subset of 
$X$ since $\OO_X^{1/q}$ is a coherent $\OO_X$-module (we always assume that 
$A$ is F-finite). Also $\tilde{\phi}$ induces an $\OO_{X'}$-linear map 
$\OO_{X'}^{1/q} \to \OO_{X'}$, which we denote by the same letter $\phi$. 
Then $\phi \colon \OO_{X'}^{1/q}\to \OO_{X'}$ gives an F-splitting of $X'$ 
since the composition map 
$\OO_{X'} \hookrightarrow \OO_{X'}^{1/q} \overset{\phi}\to\lra \OO_{X'}$ 
is the identity (cf. Lemma 2 of \cite{MS}).  

(4) Let $\phi$ be as in (3) and  fix an irreducible component $E_j$ of $E$. 
We want to examine if $\phi$ is defined at the generic point of $E_j$ or not. 
Let $\xi$ be the generic point of $E_j$ and $\eta$ be a regular parameter of 
$\OO_{X,\xi}$. We write $q - 1 = mr + i$ with $0 \le i < r$. Then by (3.3.1) 
we have 
$$\alpha(\phi) \in w^{-m}A(-iK_A - \rd{i\Delta}) 
               \subseteq b^{-1}w^{-m}H^0(X,\omega_X^{\otimes (-i)}).$$ 

(5) Now, we will show that $a_j \ge -1$ if $(A,\Delta)$ is F-pure. 
Assume, on the contrary, that $a_j < -1$. Then, since $ra_j \in \Z$, we have 
$a_j \le -1 - 1/r$, and for a fixed integer $s$, we can take $q \gg 0$ so that 
$-mra_j \ge q + s$. This means that if $\kappa$ is a local generator of the 
$\OO_X^{1/q}$-module ${\Cal H}om_{\OO_X}(\OO_X^{1/q},\OO_X)$ at $\xi$, then 
the stalk $\phi_{\xi}$ of $\phi$ at $\xi$ lies in $[\eta^q\OO_{X, \xi}]^{1/q}
\kappa$. To see this let $s = v_{E_j}(b)$, the value of $b \in A$ at $\xi$, 
and let $-mra_j \ge q + s$. Then $\OO_{X.\xi}^{1/q} \cdot \alpha_{\xi}(\kappa) 
= [w^{-m}\eta^{mra_j}\omega_{X,\xi}^{\otimes (-i)}]^{1/q} \supseteq [w^{-m}
\eta^{-q-s}\omega_{X,\xi}^{\otimes (-i)}]^{1/q}$ by (2), and 
$\alpha(\phi)_{\xi} \in [\eta^{-s}w^{-m}\omega_{X,\xi}^{\otimes (-i)}]^{1/q}$ 
by (4). It follows that $\phi_{\xi} \in [\eta^q\OO_{X, \xi}]^{1/q}\kappa$, 
whence $\phi_{\xi}(\OO_{X,\xi}^{1/q}) \subseteq \eta \OO_{X,\xi}$. This 
implies that $\xi \in X'$ but $\phi$ is not an F-splitting at $\xi$. This 
contradiction concludes that $a_i \ge -1$. 

(6) If $(A,\Delta)$ is divisorially F-regular, then for any $c \ne 0 \in 
A$ which is not in any minimal prime of $A(-\rd{\Delta})$, there exists an 
$A$-homomorphism $\psi \colon A((q-1)\Delta)^{1/q} \to A$ sending $c^{1/q}$ 
to $1$. Then $\phi = c^{1/q}\psi$ gives a splitting of $A \hookrightarrow 
A((q-1)\Delta)^{1/q}$, and we may think $\phi \in c^{1/q}\Hom_A(A((q-1)
\Delta)^{1/q},A) \cong [c \cdot A((1-q)K_A-\rd{(q-1)\Delta})]^{1/q}$. 
Let us choose $c$ so that the value $t = v_{E_j}(c)$ of $c$ at $\xi$ 
satisfies $t \ge r+s = r+v_{E_j}(b)$. Then, arguing as in the F-pure case, 
we see that if $a_j = -1$, then $\phi_{\xi} \in [\eta^q \OO_{X,\xi}]^{1/q}
\kappa$ and again $\phi$ cannot be an F-splitting at $\xi$. 

This completes the proof of the theorem. 
\enddemoo

\demo{Second proof}
Next we shall give an alternative proof of Theorem 3.3 with different flavor. 
As we have seen in the first proof, the adjunction formula gives 
$${\Cal H}om_{\OO_X}(\OO_X^{1/q},\OO_X) \cong 
               (\omega_X^{\otimes (1-q)})^{1/q} \cong \OO_X((1-q)K_X)^{1/q},$$ 
$$\Hom_A(A((q-1)\Delta)^{1/q},A) \cong A((1-q)K_A - \rd{(q-1)\Delta})^{1/q}.$$ 

(1) Now assume that $(A,\Delta)$ is F-pure and let $\phi \colon A((q-1)\Delta
)^{1/q} \to A$ be a splitting of $A \hookrightarrow A((q-1)\Delta)^{1/q}$. 
We regard $\phi \in \Hom_A(A((q-1)\Delta)^{1/q},A) \cong A((1-q)K_A - \rd{(q-1)
\Delta})^{1/q}$ 
as a rational section of the sheaf $\omega_X^{\otimes (1-q)}$ 
and consider the corresponding divisor on $X$, 
$$D = D_{\phi} = (\phi)_0 - (\phi)_{\infty},$$ 
where $(\phi)_0$ (resp.\ $(\phi)_{\infty}$) is the divisor of zeros (resp.\ 
poles) of $\phi$ as a rational section of $\omega_X^{\otimes (1-q)}$. 
Clearly, $D$ is linearly equivalent to $(1-q)K_X$ and $(\phi)_{\infty}$ 
is an $f$-exceptional divisor. Hence $f_*D$ is linearly equivalent to 
$(1-q)K_Y$ and $f_*D \ge \rd{(q-1)\Delta}$. We denote $X' = X \setminus 
\text{Supp }(\phi)_{\infty}$. Then $\phi$ lies in 
$$\Hom_{\OO_{X'}}(\OO_{X'}^{1/q},\OO_{X'}) 
                                \cong H^0(X',\omega_{X'}^{\otimes (1-q)}),$$ 
and gives an F-splitting of $X'$. 

(2) We show that the coefficient of $D$ in each irreducible component is 
$\le q-1$. Assume to the contrary that there exists an irreducible component 
of $D$, say $D_0$, whose coefficient is $\ge q$. Then $\phi$ lies in 
$$\Hom_{\OO_{X'}}(\OO_{X'}(qD_0)^{1/q},\OO_{X'}) 
                          \cong H^0(X',\omega_{X'}^{\otimes (1-q)}(-qD_0)),$$ 
and gives a splitting of the map $\OO_{X'} \hookrightarrow \OO_{X'}(qD_0)^{1/q
}$. But this map factors through $\OO_X(D_0)$, 
and $\OO_{X'} \hookrightarrow \OO_{X'}(D_0)$ never splits as an $\OO_{X'}
$-module homomorphism, since $D_0$ intersects $X'$. Consequently, every 
coefficient of $D$ must be $\le q-1$. 

(3) Let 
$$B = \frac1{q-1}D - \tilde{\Delta}.$$ 
Then $B$ is $\Q$-linearly equivalent to $-(K_X + \tilde{\Delta})$, so that 
$f_*B$ is $\Q$-linearly equivalent to $-f_*(K_X + \tilde{\Delta}) = -(K_Y + 
\Delta)$. Hence $f_*B$ is $\Q$-Cartier, and we can define the pull-back $f^*
f_*B$. Since $B + \sum_{i=1}^s a_i E_i$ is $\Q$-linearly equivalent to $-f^*
(K_Y + \Delta)$, $(B - f^*f_*B) + \sum_{i=1}^s a_i E_i$ is an $f$-exceptional 
divisor which is $\Q$-linearly trivial relative to $f$. Hence 
$$(B - f^*f_*B) + \sum_{i=1}^s a_i E_i = 0. \tag{3.3.3}$$ 

(4) Now $f_*D - (q-1)\Delta \ge \rd{(q-1)\Delta} - (q-1)\Delta \ge -\Delta'$ 
for some (effective $\Q$-) Cartier divisor $\Delta'$ on $Y$ which is 
independent of $q$. This implies $f_*B \ge -\dfrac1{q-1}\Delta'$, whence 
$f^*f_*B \ge - \dfrac1{q-1}f^*\Delta'$. Therefore, if we take $q = p^e$ 
sufficiently large, then the coefficient of $f^*f_*B$ in $E_i$ is greater 
than $-1/r$. On the other hand, we have seen in (2) that the coefficient of 
$B$ in $E_i$ is at most 1. Since $ra_i \in \Z$, it follows from (3.3.3) that 
$a_i \ge -1$. 

(5) Assume now that $(A,\Delta)$ is divisorially F-regular. For a $\Q$-Cartier 
divisor $\Delta'$ as in (4), we choose $c \ne 0 \in A$ which is not in any 
minimal prime of $A(-\rd{\Delta})$ so that the value $v_{E_i}(c)$ is greater 
than the coefficient of $f^*\Delta'$ in $E_i$. Then there exists an $A$-linear 
map $\psi\colon A((q-1)\Delta)^{1/q} \to A$ sending $c^{1/q}$ to $1$, and $\phi
= c^{1/q}\psi$ gives a splitting of $A \hookrightarrow A((q-1)\Delta)^{1/q}$. 
Let $D_{\phi}$ and $D_{\psi}$ be the divisors defined by $\phi$ and $\psi$ 
as rational sections of $\omega_X^{\otimes (1-q)}$, respectively. Then the 
coefficient of $D_{\phi}$ in $E_i$ is $\le q-1$ by (2), so that the coefficient
of $D_{\psi} + f^*\Delta'$ in $E_i$ is $< q - 1$, since $D_{\phi} = D_{\psi} + 
\text{div}_X(c)$. Let $B = (1/(q-1))D_{\psi} - \tilde{\Delta}$ and argue as in 
the F-pure case. Then we have $a_i > -1$, as required. 
\enddemoo

\remark{Remark 3.4}
The both proofs show that if the map $A \to A((q-1)\Delta)^{1/q}$ sending $1$ 
to $c^{1/q}$ has a splitting for a {\it single} element $c$ in a sufficiently 
high power of $H^0(X,\OO_X(-E))$, then we have $a_j > -1$ for every $j$. In 
the case $\Delta = 0$, if the map $A \overset{c^{1/q}}\to\lra A((q-1)\Delta
)^{1/q}$ splits for an element $c$ with $v_{E_i}(c) > 0$, then $a_i > -1$. 
\endremark

\remark{Remark 3.5}
It follows that if $(A,\Delta)$ is strongly F-regular and if $K_Y+\Delta$ is 
$\Q$-Cartier, then $\rd{\Delta} = 0$ and in (3.1.1), one has $a_j > -1$ for 
every $j$. But the above proof says something about strong F-regularity even 
when $K_Y+\Delta$ is not $\Q$-Cartier. 

Let $f \colon X \to Y = \Spec A$ be a resolution of singularities admitting 
an $f$-ample Cartier divisor $H$ supported on the exceptional locus of $f$, 
and assume that $(A,\Delta)$ is strongly F-regular. Then we can prove that 
there exists a $\Q$-divisor $G$ on $X$ satisfying the following conditions: 
\roster
\item"{(i)}" $\rup{G}$ is an effective $f$-exceptional divisor; 
\item"{(ii)}" $-f_*G \ge \Delta$; 
\item"{(iii)}" $G - K_X$ is $f$-ample. 
\endroster
To see this, for an effective $\Q$-divisor $\Delta'$ as in (4) above, choose 
$c \ne 0 \in A$ such that $\text{div}_X(c) \ge f^*\Delta'$ and that this is a 
strict inequality for the coefficients of $E_j$'s. Then there is an $A$-linear 
map $\theta \colon A((q-1)\Delta)^{1/q} \to A$ sending $c^{2/q}$ to $1$, by 
strong F-regularity. Let $\psi = c^{1/q}\theta$ and define the divisor 
$D_{\psi}$ as in the step (5) above. Now we put $G = -\dfrac1{q-1}D_{\psi} 
+ \varepsilon H$ for $\varepsilon \in \Q$ with $0 < \varepsilon \ll 1$. 
Then $G$ satisfies conditions (ii) and (iii). Also, the coefficient of $G$ 
in each irreducible component $G_j$ is $\ge -1$ and is $> -1$ if $G_j$ is 
$f$-exceptional (resp.\ $\le 0$ if $G_j$ is not $f$-exceptional). Since 
$\rd{\Delta} = 0$ and ampleness is an open condition, we can perturb 
coefficients of $G$ slightly so that $G$ satisfies condition (i) as 
well as (ii) and (iii). 
\endremark

\demo{Example 3.6}
Let us take a look at two typical examples in the case $\Delta = 0$. 

(1) Let $A$ be the localization of $k[X,Y,Z,W]/(X^4+Y^4+Z^4+W^4)$ at the 
unique graded maximal ideal, where $k$ is a field of characteristic $p$. 
By a criterion by R.~Fedder \cite{F}, $A$ is F-pure if (and only if) $p 
\equiv 1$ (mod 4).  On the other hand, if $f \colon X \to \Spec A$ is 
the blowing-up of the maximal ideal of $A$, then $X$ is regular with 
$\omega_X \cong \OO_X(-E_0)$, where $E_0$ is the exceptional divisor 
of $f$. Then the discrepancy of $E_0$ is $a_0 = -1$. 

(2) Let $A$ be the localization of the $r$th Veronese subring $k[X_1,\ldots ,
X_n]^{(r)}$ of $k[X_1,\ldots ,X_n]$ at the unique graded maximal ideal, where 
$k$ is a field of characteristic $p$.  Then $A$ is $\Q$-Gorenstein with index 
$r/(r,n)$ and is surely strongly F-regular being a pure subring of a regular 
ring.  Again, let $X$ be the blowing-up of the maximal ideal of $A$. Then $X$ 
is regular and if $E_0$ is the exceptional divisor, we have $a_0 = -1 + n/r$. 
\enddemo

Now, let us discuss rings essentially of finite type over a field of 
characteristic zero. Our goal is the following 

\proclaim{Theorem 3.7}
Let $A$ be a normal local ring essentially of finite type over a field $k$ 
of characteristic zero and let $\Delta$ be an effective $\Q$-Weil divisor on 
$Y = \Spec A$ such that $K_Y + \Delta$ is $\Q$-Carier. If $(A,\Delta)$ is of 
strongly F-regular (resp.\ divisorially F-regular, F-pure) type, that means 
the reduction modulo $p$ of $(A,\Delta)$ is strongly F-regular (resp.\ 
divisorially F-regular, F-pure) for infinitely many prime $p$, then 
$(Y,\Delta)$ is KLT (resp.\ PLT, LC). 
\endproclaim

\demo{Proof}
Let $Z$ be a normal affine variety over $k$ such that $A \cong \OO_{Z,z}$ 
for some $z \in Z$. We may assume that $\Delta$ comes from a $\Q$-divisor 
$\Delta_Z$ such that $K_Z + \Delta_Z$ is $\Q$-Cartier. Let $g \colon \tilde{X} 
\to Z$ be a log resolution of $(Z,\Delta_Z)$ (such that $g^{-1}_*\rd{\Delta_Z}$
is smooth, in order to prove PLT). Then since $Z, \Delta_Z, g, \tilde{X}$ 
are defined by finite elements of $k$, we can choose a finitely generated 
$\Z$-subalgebra $R$ of $k$, schemes $Z_R, \tilde{X}_R$ of finite type over 
$R$, a $\Q$-divisor $\Delta_R$ and an $R$-morphism $g_R \colon \tilde{X}_R 
\to Z_R$, which give $Z, \tilde{X}, \Delta_Z$ and $g$ after tensoring $k$ 
over $R$. Then, taking suitable open subset of $\Spec R$, we may assume 
$Z_R, \tilde{X}_R$ and each irreducible component of $\Delta_R$ are flat 
over $R$, $\tilde{X}_R$ is smooth over $R$ with each irreducible component 
of the exceptional divisor and $Z_R$ is normal over $R$ (cf. \cite{EGA, IV, 
12.1.7}). Also, we can choose $R$ to be regular with trivial dualizing module. 
Then we may assume that the numbers $a_i$ in the formula  
$K_X + \tilde{\Delta} = f^*(K_Y + \Delta) + \sum _{i=1}^r a_iE_i$ 
are preserved in the fibers of $g_R \colon \tilde{X}_R \to Z_R$ over $\Spec R$ 
in an open neighborhood of the generic point of $\Spec R$. By our assumption, 
this open neighborhood contains a maximal ideal $\p$ of $R$, such that the 
base change of $(Z_R,\Delta_R)$ to $k(\p)$ over $R$ is strongly F-regular 
(resp.\ divisorially F-regular, F-pure) at $\bar{z}$, a specialization of 
$z$ (to be F-regular or F-pure is preserved under base field extension or 
restriction and so depends only on characteristic). This implies $a_i>-1$ 
(resp.\ $a_i \ge -1$) by Theorem 3.3 and our assertion is proved.  
\enddemoo

\head 4. Applications \endhead

\demo{4.1. The graded case} (cf.\ \cite{W1, W2})
Let $R = \bigoplus_{n\ge 0} R_n$ be a normal graded ring over a perfect field 
$R_0 = k$ of characteristic $p > 0$. Given a homogeneous element $T$ of degree 
$1$ in the quotient field of $R$, there is an ample $\Q$-Cartier divisor $D$ 
on $X = \Proj R$ such that 
$$R = R(X,D) = \bigoplus_{n\ge 0} H^0(X,\OO_X(nD))T^n.$$ 
If $D = \sum_{i=1}^s (e_i/d_i)D_i$ for distict prime divisors $D_1,\dots ,D_s$ 
and coprime integers $d_i$ and $e_i$ with $d_i > 0$ ($1 \le i \le s$), we 
denote $D'= \sum_{i=1}^s ((d_i-1)/d_i)D_i$ and call it the ``fractional part" 
of $D$ (\cite{W1}). Then $K_R^{(r)}$ is free if and only if $r(K_X+D')$ is 
linearly equivalent to $bD$ for some integer $b$. 

If $Z = \Spec_X \left(\bigoplus_{n\ge 0} \OO_X(nD)T^n\right)$ is the 
``graded blowing-up" of $\Spec R$ and if $E_0 \cong X$ is the exceptional 
divisor of this blowing-up, its discrepancy is $a_0 = -1 - b/r$. 
Consequently, if $R$ is strongly F-regular (resp.\ F-pure) then $b < 0$ 
(resp.\ $b \le 0$). (In this case, although $Z$ is not Gorenstein in general, 
the number $a_0$ is preserved after we make more blowing-ups and get a 
Gorenstein (or regular) scheme.) 

Now let $\Delta = \sum_{j=1}^r t_j\Delta_j$ be an effective $\Q$-Weil divisor 
on $\Spec R$ which is stable under the $k^*$-action, i.e., each irreducible 
component $\Delta_j$ of $\Delta$ is defined by a {\it homogeneous} prime ideal 
of $R$ of height $1$. Then there exist effective $\Q$-Weil divisors $\Gamma_1,
\dots, \Gamma_r$ on $X = \Spec R$ such that for every $i \in \Z$ and $j = 1,
\dots ,r$, 
$$R(i\Delta_j) = \bigoplus_{n\in \Z} H^0(X,\OO_X(i\Gamma_j + nD))T^n.$$ 
Let $\Gamma = \sum_{j=1}^r t_j \Gamma_j$ and $D'$ be the ``fractional part" 
of $D$ as above. Then we can rephrase Proposition 2.4 using this language, 
cf.\ \cite{W2}. 
\enddemo

\proclaim{Proposition 4.2}
Let the notation be as in 4.1 and let $\dim R = d+1 \ge 2$. 
\roster
\item $(R,\Delta)$ is F-pure if and only if the map 
$$F^e \colon H^d(X,\OO_X(K_X)) \to H^d(X,\OO_X(q(K_X+D')+(q-1)\Gamma))$$ 
is injective for every $q = p^e$. 
\item $(R,\Delta)$ is divisorially $($resp.\ strongly$)$ F-regular if and 
only if for every $n\ge 0$ and every $c \ne 0 \in H^0(X,\OO_X(nD)) \setminus 
\bigcup_{t_j\ge 1} H^0(X,\OO_X(-\Gamma_j+nD))$ $($resp.\ $c\ne 0 \in H^0(X,
\OO_X(nD)))$, there exists $q = p^e$ such that the map 
$$\phantom{mmm}
  cF^e \colon H^d(X,\OO_X(K_X)) \to H^d(X,\OO_X(q(K_X+D')+(q-1)\Gamma + nD))$$ 
is injective. 
\endroster
\endproclaim

\demo{Example 4.3}
Let $R = R(\PP^1,D)$ with $D = \frac12 ((0)+(\infty))$ and let $\Delta$ be the 
divisor defined by $R(-\Delta)= \bigoplus_{n\ge 0} H^0(\PP^1,\OO(-(1)+nD))T^n$.
Then $(R,\Delta)$ is F-pure if and only if $p \ne 2$. This is a special case 
of the type (c) in 4.4 below. 
\enddemo

\demo{4.4.\ Two-dimensional F-regular and F-pure pairs with 
integer coefficient boundary}\hfill

\noindent
Combining Theorem 3.3 and the technique used in \cite{Ha2}, we can classify 
F-regular and F-pure pairs $(A,\Delta)$ such that $A$ is a two-dimensional 
normal local ring of characteristic $p > 0$ and that $\Delta$ is an integer 
coefficient effective Weil divisor on $\Spec A$. The case $\Delta = 0$ is 
treated in \cite{Ha2}, so we treat here only the case $\Delta \ne 0$. 

Two-dimensional LC pairs with nonzero reduced boundaries are classified in 
terms of the dual graph of the union of the exceptional divisor $E$ and the 
strict transform $\tilde{\Delta}$ of $\Delta$ for the minimal log resolution 
$f \colon X \to Y = \Spec A$, and such a pair is one of the following three 
types (\cite{A}, \cite{Ka}): 

\newpage
\roster
\item"{(a)}" $\bullet - \circ - \cdots - \circ$
\item"{(b)}" $\bullet - \circ - \cdots - \circ - \bullet$
\item"{(c)}" 
$\matrix
\phantom{\bullet \! - \circ - \cdots \! -}\circ -2 \phantom{- \circ} \\
\phantom{\bullet \! - \circ - \cdots -}\vert \phantom{- \circ -\! 2} \\
         \bullet - \circ - \cdots - \circ - \circ -\! 2
\endmatrix $ 
\endroster
Here, a blank circle $\circ$ (resp.\ a solid circle $\bullet$) denotes an 
irreducible component $E_i \cong \PP^1$ of $E$ (resp.\ an irreducible component
of $\tilde{\Delta} = f^{-1}_*(\Delta)$), and the numbers ``$-2$" outside the 
circles in (c) mean that the corresponding components have self-intersection 
number equal to $-2$. Among the above three types, type (a) is PLT but types 
(b) and (c) are not PLT. 

Theorem 3.3 tells us that a two-dimensional F-pure (resp.\ divisorially 
F-regular) pair with reduced boundary is of type (a), (b) or (c) (resp.\ type 
(a)) as above. We can show that the converse is true but one exception. 
\enddemo

\proclaim{Theorem 4.5}
Let $(A,\m)$ be a two-dimensional normal local ring essentially of finite type 
over an algebraically closed field $k = A/\m$ of characteristic $p > 0$ and 
let $\Delta$ be a nonzero reduced Weil divisor on $Y = \Spec A$. Let $f \colon 
X \to Y$ be the minimal resolution with exceptional divisor $E$ and let 
$\tilde{\Delta} = f^{-1}_*(\Delta)$. 
\roster
\item $(A,\Delta)$ is divisorially F-regular if and only if the dual graph of 
$E \cup \tilde{\Delta}$ is of type {\rm (a)} in 4.4. 
\item $(A,\Delta)$ is F-pure but not divisorially F-regular if and only if the 
dual graph of $E \cup \tilde{\Delta}$ is of type {\rm (b)} in 4.4, or of type 
{\rm (c)} in 4.4 and $p \ne 2$. 
\endroster
\endproclaim

\demo{Proof}
(1) Let the dual graph of $E \cup \tilde{\Delta}$ be of type (a) as in 4.4, 
and fix any nonzero element $c \in A$ which is not in $A(-\Delta)= H^0(X,\OO_X
(-\tilde{\Delta}))$. By Proposition 2.4, our goal is to show that there exists 
$q = p^e$ such that the map 
$$cF^e \colon H_{\m}^2(K_A) \to H_{\m}^2(A(qK_A+(q-1)\Delta)) \tag{4.5.1}$$ 
is injective. Blowing up at $E \cap \tilde{\Delta}$ finitely many times if 
necessary, we may assume that $\tilde{\Delta}$ does not intersect the strict 
transform of $\Div_Y(c)$ on $X$ and that the graph of $E \cup \tilde{\Delta}$ 
is still of type (a) as follows. 
$$\tilde{\Delta}\text{---}E_0\text{---}E_1\text{---}\cdots\text{---}E_l$$ 
Note that $f \colon X \to Y$ may not be the minimal resolution any longer, but 
$E_i$ is a $(-1)$-curve only if $i = 0$. By Lemma 3.9 of \cite{Ha2}, one has 
an effective $f$-exceptional $\Q$-divisor 
$$D = E_0 + \frac{d-1}d E_1 + (\text{terms of }E_2, \dots ,E_l)$$ 
with $d \in \N$ such that $(K_X + D)E_0 = -1 - 1/d$ and $(K_X + D)E_j = 0$ 
for $j = 1, \dots ,l$. Let $D' = D - E_0$. 

Let $n$ be the integer with $c \in H^0(X,\OO_X(-nE_0)) \setminus H^0(X,\OO_X
(-(n+1)E_0))$. Since $H_{\m}^2(K_A) \cong H_E^2(\omega_X)$ by $H^i(X,\omega_X) 
= 0$ $(i = 1,2)$, the map (4.5.1), followed by 
$H_{\m}^2(A(qK_A+(q-1)\Delta)) \to 
                 H_E^2(\OO_X(q(K_X+D')+(q-1)(\tilde{\Delta}+E_0)-nE_0))$, is 
$$cF^e \colon H_E^2(\omega_X) \to 
       H_E^2(\OO_X(q(K_X+D')+(q-1)(\tilde{\Delta}+E_0)-nE_0)), \tag{4.5.2}$$ 
and it is sufficient to show that this map is injective for some $q = p^e$. 

Let $P,Q \in E_0$ be the points of intersection of $E_0$ with $E_1,\tilde{
\Delta}$, respectively, and let $\frak{d}' = D'|_{E_0} = \frac{d-1}d P$ and 
$\frak{a} = -E_0|_{E_0}$ as ($\Q$-)divisors on $E_0 \cong \PP^1$. Then $c \in 
H^0(X,\OO_X(-nE_0))$ restricts to a nonzero element $\bar{c} \in H^0(E_0,
\OO_{E_0}(n\frak{a}))$, and $\bar{c} \notin H^0(E_0,\OO_{E_0}(n\frak{a}-Q))$ 
since $\tilde{\Delta} \cap f^{-1}_*\Div_Y(c) = \emptyset$. 

Now if $q = p^e$ is 
sufficiently large, then $-q(K_X+D)-(q-1)\tilde{\Delta}+nE_0$ is $f$-nef, so 
that $R^1f_*\omega_X(\rup{-q(K_X+D)-(q-1)\tilde{\Delta}+nE_0}) = 0$ (cf.\ 
\cite{Ha2, Lemma 3.3}), or dually, $H_E^1(\OO_X(q(K_X+D)+(q-1)\tilde{\Delta}-
nE_0)) = 0$. Hence we have the following commutative diagram with the vertical 
arrows being injective for $q = p^e \gg 0$. 
$$\CD
H_E^2(\omega_X) @>cF^e>> H_E^2(\OO_X(q(K_X+D')+(q-1)(\tilde{\Delta}+E_0)-nE_0))
\\ @AAA @AAA \\
H^1(E_0,\omega_{E_0}) @>\bar{c}F^e>> 
               H^1(E_0,\OO_{E_0}(q(K_{E_0} + \frak{d}') + (q-1)Q + n\frak{a})) 
\endCD
$$
Note also that the map $H^1(E_0,\omega_{E_0}) \to H_E^2(\omega_X)$ on the left 
is identified with $k = A/\m \hookrightarrow E_A(A/\m)$, whence an essential 
extension. By computing \v{C}ech cohomologies on $E_0 \cong \PP^1$, we can 
verify that the map $\bar{c}F^e$ at the bottom is injective for $q = p^e \gg 
0$, from which follows the injectivity of the map (4.5.2) at the top. 

(2) If the graph is of type (b) or (c), we use the anti-discrepancy of $f$ in 
place of $D$ in the proof of (1) above, i.e., $D = f^*(K_Y + \Delta) - (K_X +
\tilde{\Delta})$. Let $D_0= \rd{D}$, $D'= D- D_0$ and $\frak{d}'= D'|_{D_0}$, 
and argue as in (1) for $c = 1$. Then it follows that $(A,\Delta)$ is F-pure 
if and only if the map 
$$F^e \colon H^1(D_0,\omega_{D_0}) \to 
    H^1(D_0,\OO_{D_0}(q(K_{D_0} + \frak{d}') + (q-1)\tilde{\Delta}|_{D_0}))$$ 
is injective for all $q = p^e$. We can show as in \cite{Ha2, Claim 4.8.1} that 
this map is injective if and only if the graph of $E \cup \tilde{\Delta}$ is 
of type (b), or of type (c) and $p \ne 2$. See also \cite{Ha2} for details. 
\enddemoo

\demo{4.6.\ Strongly F-regular rings vs.\ admissible singularities}
The correspondence of $\Q$-Gorenstein (strongly) F-regular rings and log 
terminal singularities (in the case $\Delta = 0$) is now well established. 
Namely, we have seen that a ring of characteristic zero has log terminal 
singularities if it is of (strongly) F-regular type and $\Q$-Gorenstein. 
The converse of this implication is also proved to be true \cite{Ha3}. 
Thus we are tempted to consider what non-$\Q$-Gorenstein strongly F-regular 
rings are. We have no established answer to this question, but there is a 
candidate which is expected to correspond strongly F-regular rings even in 
non-$\Q$-Gorenstein case. 

In \cite{N}, Nakayama introduced the notion of admissible singularities as 
an analog of log terminal singularities in the absence of $\Q$-Gorensteinness. 
Let $Y$ be a normal variety of characteristic zero and let $\Delta$ be an 
effective $\Q$-Weil divisor on $Y$. Then the pair $(Y,\Delta)$ is said to be 
{\it strictly admissible} if there exist a birational morphism $f \colon X 
\to Y$ from a nonsingular $X$ and a $\Q$-divisor $G$ on $X$ satisfying the 
following conditions: 

\newpage
\roster
\item"{(i)}" $\rup{G}$ is an effective $f$-exceptional divisor; 
\item"{(ii)}" $-f_*G \ge \Delta$; 
\item"{(iii)}" $G - K_X$ is $f$-ample; 
\item"{(iv)}" $\text{Supp}(G-\rd{G})$ is a normal crossing divisor. 
\endroster

On the other hand, we observed in Remark 3.5 that if $(Y,\Delta)$ is a 
strongly F-regular pair in characteristic $p > 0$ and $f \colon X \to Y$ is 
a resolution, then we can construct a $\Q$-divisor $G$ satisfying conditions 
(i), (ii) and (iii). But we do not have condition (iv), and even worse, this 
construction depends on characteristic $p$. However, the following example 
may be a positive evidence to the correspondence of strong F-regularity and 
admissible singularity. 
\enddemo

\demo{Example 4.7}
Let $X$ be a smooth Fano variety of characteristic zero (i.e., a smooth 
projective variety with ample anti-canonical divisor $-K_X$), and $D$ be 
any ample Cartier divisor on $X$. Let $R = R(X,D)$ and let $f \colon Z \to 
\Spec R$ be the ``graded blowing-up" as in (4.1). 
The exceptional set of $f$ is a smooth divisor $E \cong X$, so if we put $G 
= (\varepsilon -1)E$ with $0 < \varepsilon \le 1$, conditions (i), (ii), (iv) 
in (4.6) are satisfied. Also, if we choose $\varepsilon \ll 1$, then $(G - K_Z)
\vert_E = (-(K_Z+E) + \varepsilon E)\vert_E = -K_E - \varepsilon D$ is ample 
since $E \cong X$ is Fano, so that condition (iii) is satisfied. Consequently, 
$R$ has an admissible singularity whereas it is not $\Q$-Gorenstein in general.

On the other hand, $R(X,-K_X)$ has a Gorenstein log terminal singularity, 
so that it is of F-regular type \cite{Ha3}. However, since the strong 
F-regularity of $R(X,D)$ depends only on the ``fractional part" $D'$ of 
$D$ \cite{W2} and the fractional parts of $D$ and $-K_X$ are both $0$ in 
this case, $R = R(X,D)$ also has strongly F-regular type. 
\enddemo

We study more about similarity of F-regularity and F-purity of pairs and 
``singularities of pairs" in characteristic zero \cite{Ko}, \cite{Sh}. The 
following theorem generalizes \cite{W2, Theorem 2.7}. See Remark 1.2 (5). 

\proclaim{Theorem 4.8}
Let $(A,\m) \to (B,\frak{n})$ be a finite local homomorphism of F-finite 
normal local rings which is \'etale in codimension 1. Let $\Delta_A$ be an 
effective $\Q$-Weil divisor on $\Spec A$ and let $\Delta_B = \pi^*\Delta_A$ 
be the pull-back of $\Delta_A$ by the induced morphism $\pi \colon \Spec B 
\to \Spec A$. If $(A,\Delta_A)$ is F-pure $($resp.\ divisorially or strongly 
F-regular$)$, then so is $(B,\Delta_B)$, too. 
\endproclaim

\demo{Proof}
Let $d = \dim A = \dim B$. Since $A \to B$ is \'etale in codimension 1, the 
natural maps $A^{1/q} \otimes_A B \to B^{1/q}$ and $A(qK_A + (q-1)\Delta_A) 
\otimes_A B \to B(qK_B+(q-1)\Delta_B)$ are isomorphic in codimension 1. Hence, 
by \cite{W2, Lemma 2.2}, we have 
$$\aligned
H_{\m}^d(A(qK_A+(q-1)\Delta_A)^{1/q}) & \otimes_A B \\
 & \cong H_{\m}^d(A(qK_A+(q-1)\Delta_A)^{1/q} \otimes_A B) \\ 
 & \cong 
    H_{\m}^d(A(qK_A+(q-1)\Delta_A)^{1/q}\otimes_{A^{1/q}}A^{1/q}\otimes_A B) \\
 & \cong H_{\m}^d(A(qK_A+(q-1)\Delta_A)^{1/q} \otimes_{A^{1/q}} B^{1/q}) \\
 & \cong H_{\m}^d((A(qK_A+(q-1)\Delta_A) \otimes_A B)^{1/q}) \\
 & \cong H_{\frak{n}}^d(B(qK_B+(q-1)\Delta_B)^{1/q}). 
\endaligned
$$
Therefore, tensoring the map 
$$F^e \colon H_{\m}^d(K_A) \to H_{\m}^d(A(qK_A+(q-1)\Delta_A)) \tag{4.8.1}$$ 
with $B$ over $A$ yields 
$$F^e \colon H_{\frak{n}}^d(K_B) \to H_{\frak{n}}^d(B(qK_B+(q-1)\Delta_B)). 
                                                               \tag{4.8.2}$$ 

Now, if $(A,\Delta_A)$ is F-pure, then for every $q = p^e$, the map (4.8.1) 
is injective, and even a {\it splitting} injective map, since $H_{\m}^d(K_A) 
\cong E_A(A/\m)$ is an injective $A$-module. Hence the map (4.8.2) is also 
injective for every $q = p^e$, so that $(B,\Delta_B)$ is F-pure. 

Similarly, if $(A,\Delta_A)$ is strongly F-regular, then for every $c \ne 0 
\in A$, there exists $q = p^e$ such that the map 
$$cF^e \colon H_{\frak{n}}^d(K_B) \to H_{\frak{n}}^d(B(qK_B+(q-1)\Delta_B))$$ 
is injective. This is true for every $c \ne 0 \in B$, since $A \to B$ is a 
finite extension of normal domains, so that $cB \cap A \ne 0$. Hence $(B,
\Delta_B)$ is strongly F-regular. 

To prove the assertion for divisorial F-regularity, note     %that $\pi$ is 
%unramified, in particular, $\pi^*\rd{\Delta_A} = \rd{\Delta_B}$. This implies 
that if $c \in B$ is not in any minimal prime ideal of $B(-\rd{\Delta_B})$, 
then there is an element of $cB \cap A$ which is not in any minimal prime 
ideal of $A(-\rd{\Delta_A})$. Then the argument for strong F-regularity 
works also for divisorial F-regularity. 
\enddemoo

\proclaim{Theorem 4.9} {\rm (cf.\ \cite{Ko, Theorem 7.5})}
Let $(A,\m)$ be an F-finite normal local ring of characteristic $p > 0$ with 
$Y = \Spec A$ and let $x \in \m$ be a nonzero element. 
\roster
\item If the pair $(A,\Div_Y(x))$ is divisorially F-regular, then the ring 
$A/xA$ is strongly F-regular. 
\item Assume that $A$ is $\Q$-Gorenstein and that the order $r$ of the 
canonical class in the divisor class group $\operatorname{Cl}(A)$ is not 
divisible by $p$. Then, if the ring $A/xA$ is strongly F-regular, the pair 
$(A,\Div_Y(x))$ is divisorially F-regular. 
\endroster
\endproclaim

\demo{Proof}
Let $B = A/xA$ and $d = \dim A$. We will look at the Frobenius actions on 
$E_A = E_A(A/\m) \cong H_{\m}^d(K_A)$ and $E_B = E_B(B/\m B) \cong H_{\m}^{d-1}
(K_B)$ in slightly different ways in proving (1) and (2), respectively. 

(1) Assume that $(A,\Div_Y(x))$ is divisorially F-regular. Then for every 
element $c \in A$ which is not in any minimal prime ideal of $xA$, there 
exists $q = p^e$ such that the map $cF^e \colon H_{\m}^d(K_A) \to H_{\m}^d
(A(qK_Y+(q-1)\Div_Y(x)))$ is injective, or equivalently, the map 
$$\alpha \colon E_A = E_A \otimes_A A \to E_A \otimes_A A^{1/q}$$ 
sending $\xi \in E_A$ to $\xi\otimes (cx^{q-1})^{1/q} \in E_A\otimes_A A^{1/q}$
is injective (see the proof of Proposition 2.6). On the other hand, we have 
an inclusion $\imath \colon E_B \cong (0:x)_{E_A} \hookrightarrow 
E_A$, and this map gives rise to a well-defined $A^{1/q}$-homomorphism 
$$\jmath \colon E_B \otimes_B B^{1/q} \to E_A \otimes_A A^{1/q}$$
sending $\xi \otimes \bar{a}^{1/q} \in E_B \otimes_B B^{1/q}$ to $\imath(\xi) 
\otimes (ax^{q-1})^{1/q} \in E_A \otimes_A A^{1/q}$, where $\bar{a}$ denotes 
the image of $a \in A$ in $B = A/xA$. Then we have the following commutative 
diagram. 
$$\CD
            E_B             @>\imath>>          E_A            \\
@V 1\otimes\bar{c}^{1/q} VV                 @V \alpha VV       \\
   E_B \otimes_B B^{1/q}    @>\jmath>>  E_A \otimes_A A^{1/q}
\endCD
$$
Hence the injectivity of the map $\alpha$ implies that 
$1_{E_B}\otimes\bar{c}^{1/q} \colon E_B \to E_B \otimes_B B^{1/q}$ is 
injective. Consequently, we have that $B$ is strongly F-regular by Remark 
1.2 (1). 

(2) Assume that $B$ is strongly F-regular. Then $B$ is a Cohen--Macaulay 
normal domain, so that $A$ and $K_A$ are also Cohen--Macaulay. 
For any $q = p^e$ and $c \in A \setminus xA$, we consider the following 
commutative diagram with exact rows, 
$$\CD
0 @>>>    K_A    @>x>>    K_A    @>>>      K_A/xK_A        @>>> 0 \phantom{,}\\
& & @Vcx^{q-1}F^eVV     @VcF^eVV           @VcF^eVV       &                  \\
0 @>>> K_A^{(q)} @>x>> K_A^{(q)} @>>> K_A^{(q)}/xK_A^{(q)} @>>> 0,
\endCD
$$
where $K_A/xK_A \cong K_B$. Since $H_{\m}^{d-1}(K_A^{(q)}/xK_A^{(q)}) \cong 
H_{\m}^{d-1}(K_B^{(q)})$ by the normality of $B$, this diagram induces the 
following commutative diagram of local cohomologies. 
$$\CD
   H_{\m}^{d-1}(K_B)    @>>>    H_{\m}^d(K_A)    \\ 
   @V cF^e VV                @Vcx^{q-1}F^eVV     \\
H_{\m}^{d-1}(K_B^{(q)}) @>>> H_{\m}^d(K_A^{(q)})
\endCD
$$
Here, the map $H_{\m}^{d-1}(K_B) \to H_{\m}^d(K_A)$ upstairs is an essential 
extension, and if $q \equiv 1$ (mod $r$), then the map $H_{\m}^{d-1}(K_B^{(q)})
\to H_{\m}^d(K_A^{(q)})$ at the bottom is also injective since $K_A \cong 
K_A^{(q)}$ is Cohen--Macaulay. 

Now, since $B$ is strongly F-regular, the map $cF^e \colon H_{\m}^{d-1}(K_B) 
\to H_{\m}^{d-1}(K_B^{(q)})$ in the above diagram is injective for all $q= p^e 
\gg 0$, by (2) of Propositions 2.2 and 2.4. Since $r$ is not divisible by $p$, 
we can choose $q = p^e$ so that $q \equiv 1$ (mod $r$) and that $cF^e\colon 
H_{\m}^{d-1}(K_B)\to H_{\m}^{d-1}(K_B^{(q)})$ is injective. Then the diagram 
implies that the map $cx^{q-1}F^e \colon H_{\m}^d(K_A) \to H_{\m}^d(K_A^{(q)})$
is injective for such a $q = p^e$, or equivalently, 
$$cF^e \colon H_{\m}^d(K_A) \to H_{\m}^d(A(qK_Y+(q-1)\Div_Y(x)))$$ 
is injective. Hence $(A,\Div_Y(x))$ is divisorially F-regular. 
\enddemoo

\remark{Remark 4.10}
(1) In the situation of Theorem 4.9 (2), assume in addition that $B = A/xA$ 
is normal and Cohen--Macaulay. Then we can also prove that $(A,\Div_Y(x))$ is 
F-pure if and only if $B = A/xA$ is F-pure. 

(2) The proof of Theorem 4.9 suggests a more general statement as follows: Let 
$(A,\m)$ and $x \in \m$ be as in Theorem 4.9. Let $S = \Div_Y(x) \cong \Spec 
A/xA$ and $\Delta$ be an effective $\Q$-Weil divisor on $Y = \Spec A$ such 
that $r(K_Y + \Delta)$ is Cartier for some positive integer $r$ which is not 
divisible by $p$. Then $(A,S+\Delta)$ is divisorially F-regular if and only if 
$(A/xA,\Delta |_S)$ is strongly F-regular. 

If we replace ``divisorially F-regular" and ``strongly F-regular" in this 
assertion by ``PLT" and ``KLT" respectively, we find the so-called ``inversion 
of adjunction" in characteristic zero \cite{Ko, Theorem 7.5}, \cite{Sh}. 
\endremark

\proclaim{Corollary 4.11} {\rm (cf.~\cite{AKM})}
Let $(A,\m)$ be an F-finite local ring of characteristic $p > 0$ and let 
$x \in \m$ be a nonzero element. Assume that $A$ is $\Q$-Gorenstein and 
that the order $r$ of the canonical class in the divisor class group of 
$Y = \Spec A$ is not divisible by $p$. If $A/xA$ is strongly F-regular, 
then $A$ is also strongly F-regular. 
\endproclaim

\demo{Proof}
If $A/xA$ is strongly F-regular, then $A/xA$ is F-rational. This implies that 
$A$ is F-rational, whence normal \cite{FW}, \cite{HH3, 4.2}. Thus we can apply 
Theorem 4.9, which implies that $(A,\Div_Y(x))$ is divisorially F-regular. It 
follows that $(A,0)$ is divisorially F-regular, meaning that $A$ is strongly 
F-regular. 
\enddemoo

In fact, we do not have to assume that the index $r$ of $A$ is not divisible 
by $p$ in Corollary 4.11 (see Aberbach, Katzman and MacCrimmon \cite{AKM}). 
On the other hand, an example by Singh \cite{Si} shows that Corollary 4.11 
fails in the absense of $\Q$-Gorensteinness. 

\head 5. Open Problems \endhead

So far, we have seen many positive evidences to the correspondence of 
F-purity (resp.\ divisorial, strong F-regularity) and LC (resp.\ PLT, KLT) 
property, which enable us to ask about the converse of Theorem 3.7. 

\demo{5.1}
Let $A$ be a normal ring essentially of finite type over a field of 
characteristic zero and let $\Delta$ be a $\Q$-Weil divisor on $Y = \Spec A$. 
\enddemo

\demo{Conjecture 5.1.1}
If $(Y,\Delta)$ is KLT (resp.\ PLT), then it is of open strongly F-regular 
type (resp.\ open divisorially F-regular type), i.e., the modulo $p$ reduction 
of $(A,\Delta)$ is strongly F-regular (resp.\ divisorially F-regular) for all 
$p \gg 0$. 
\enddemo

Conjecture 5.1.1 is true when $\Delta = 0$ \cite{Ha3}. Also, Theorem 4.8 and 
\cite{Ko, Theorem 7.5} shows that a pair $(A,\Div_Y(x))$ is PLT if and only if 
it is of divisorially F-regular type and $A$ is $\Q$-Gorenstein, since $A/xA$ 
is log terminal if and only if it is of F-regular type and $\Q$-Gorenstein. 

\demo{Problem 5.1.2}
If $(Y,\Delta)$ is LC, then is it always of (dense) F-pure type, i.e., the 
modulo $p$ reduction of $(A,\Delta)$ is F-pure for infinitely many $p$? 
\enddemo

\demo{5.2. Log canonical thresholds} (\cite{Ko, \S 8}, \cite{Sh})
An affirmative answer to (5.1) suggests a Frobenius characterization of an 
important invariant called the log canonical threshold. Let $Y$ be a variety 
in characteristic zero with only log terminal singularity at a point $y \in 
Y$ and $\Delta$ be an effective $\Q$-Cartier divisor on $Y$. The {\it log 
canonical threshold} of $\Delta$ at $y \in Y$ is defined by 
$$LCTh_y(Y,\Delta) = \text{sup}\{ t \in \R \; \vert \; (Y,t\Delta) 
                                             \text{ is LC at } y \in Y \}.$$ 
Theorem 3.3 tells us that this invariant is greater than or equal to the 
supremum of $t \in \R$ such that reduction modulo $p$ of $(\OO_{Y,y},t\Delta)$ 
is  F-pure for infinitely many $p$. 

\demo{Conjecture 5.2.1} In the notation as above, the following are equal 
to each other: 
\roster
\item"{(i)}" $LCTh_y(Y,\Delta)$; 
\item"{(ii)}" sup$\{ t \in \R$ $\vert$ reduction mod $p$ of $(\OO_{Y,y},
                              t\Delta)$ is F-pure for infinitely many $p$\}; 
\item"{(iii)}" sup$\{ t \in \R$ $\vert$ reduction mod $p$ of $(\OO_{Y,y},
                            t\Delta)$ is strongly F-regular for $p \gg 0$\}. 
\endroster
\enddemo

We expect that $LCTh_y(Y,\Delta)$ is computable in terms of characteristic $p$ 
method. 
For example, by a Fedder-type criterion we can show that $(k[[x,y]],t \cdot
\Div(x^2-y^3))$ is F-pure if and only if $t \le 5/6$ (Example 2.9 (3)), 
and $5/6$ is nothing but the log canonical threshold of the hypersurface 
$x^2 - y^3 = 0$ in $\Bbb{A}_k^2$ at the origin. Also, compare Proposition 2.10 
with \cite{Ko, Lemma 8.10}. 
\enddemo

\demo{5.3} Finally, we propose other open problems. 

\demo{Problem 5.3.1} If $(A,\Delta)$ is of strongly F-regular type ($K_Y + 
\Delta$ is not necessarily $\Q$-Cartier), then does the pair $(Y,\Delta)$ have 
admissible singularities? How about the converse implication? 
\enddemo

\demo{Problem 5.3.2} Let $(A,\m)$ be a $d$-dimensional normal local ring of 
characteristic $p > 0$ and let $\Delta$ be an effective $\Q$-Weil divisor on 
$Y = \Spec A$ such that $\rd{\Delta} = 0$. Define the {\it tight closure of 
zero submodule in $E = E_A(A/\m) \cong H_{\m}^d(K_A)$ with respect to the 
pair $(A,\Delta)$},  $0_E^{*\Delta} \subset E = H_{\m}^d(K_A)$, by 
$$\aligned
z \in 0_E^{*\Delta} & \Longleftrightarrow \\ 
& \exists c \ne 0 \in A \text{ such that 
          $cz^q = 0$ in $H_{\m}^d(A(qK_A+(q-1)\Delta))$ for } q= p^e \gg 0,
\endaligned
$$ 
where $z^q$ denotes the image of $z \in E$ by $F^e \colon H_{\m}^d(K_A) \to 
H_{\m}^d(A(qK_A+(q-1)\Delta))$. 
We ask for a geometric interpretation of the ideal 
$\tau(A,\Delta) = \text{Ann}_A(0_E^{*\Delta})$. 
Especially, if $K_Y + \Delta$ is $\Q$-Cartier, is $\tau(A,\Delta)$ equal to 
the multiplier ideal of the pair $(Y,\Delta)$ in characteristic $p \gg 0$? 
The answer to this question is affirmative if $\Delta = 0$ (\cite{Ha4}, 
\cite{S2}). See also \cite{HH1,3} for tight closure, and \cite{E} for 
multiplier ideals. 
\enddemo

\vskip 10pt
Department of Mathematical Sciences, Waseda University, 
Okubo, Shinjuku-ku, Tokyo 169--8555, Japan

E-mail: nhara\@mn.waseda.ac.jp

\vskip 5pt
Department of Mathematics, College of Humanities and Sciences, Nihon 
University, 
Sakura-josui, Setagaya-ku, Tokyo 156--0045, Japan

E-mail: watanabe\@math.chs.nihon-u.ac.jp

\end